\documentclass[reqno,11pt]{amsart}
\usepackage[top=2cm,bottom=2cm,left=2cm,right=2cm]{geometry}
\usepackage{amsthm,amsmath,amssymb,dsfont}
\usepackage{mathrsfs,amsfonts,functan,extarrows,mathtools}
\usepackage{cases,latexsym,mathrsfs}
\usepackage{amsmath}
\usepackage{indentfirst}
\usepackage{appendix}
\usepackage{xcolor}
\usepackage{indentfirst, latexsym, amssymb, enumerate,amsmath,graphicx}
\usepackage{float}
\usepackage{stmaryrd}
 \usepackage[colorlinks,linkcolor=black,citecolor=black]{hyperref}
\allowdisplaybreaks[4]

\numberwithin{equation}{section}
\newtheorem{theorem}{Theorem}[section]

\newtheorem{lemma}{Lemma}[section]
\newtheorem{remark}{Remark}[section]
\newtheorem{prop}{Proposition}[section]
\newtheorem{corol}{Corollary}[section]

\newcommand{\e}{{\epsilon}}

\begin{document}


\def\esssup{\mathop{\rm ess\, sup}}

\title[Combined limits of the NSF-P1 model]
{The combined non-equilibrium diffusion and low Mach number limits of the compressible Navier-Stokes-Fourier-P1 approximation radiation model}

\author[F. Li]{Fucai Li}
\address{School of Mathematics, Nanjing University, Nanjing
 210093, P. R. China}
\email{fli@nju.edu.cn}
\author[S. Zhang]{Shuxing Zhang$^\dag$}\thanks{$^\dag$ Corresponding author. 
}
\address{School of Mathematical Sciences, Jiangsu University, Zhenjiang 212013, P. R. China}
\email{zhangsx@ujs.edu.cn}

\date{}

\begin{abstract}
In this paper, we investigate the combined non-equilibrium diffusion and low Mach number limits of the compressible Navier-Stokes-Fourier-P1 (NSF-P1) model with general initial data, which arises in the radiation hydrodynamics.
Compared to the classical compressible Navier-Stokes-Fourier system, the NSF-P1 model has an asymmetric singular structure caused by the radiation field.
To handle these singular terms, we introduce an equivalent pressure and an equivalent velocity to balance the order of singularity and establish the uniform estimates of solutions by designating appropriate weighted norms as well as carrying out delicate energy analysis.
We conclude that, for partially general initial data and the strong scattering effect, the NSF-P1 model converges to the system of low Mach number heat-conducting viscous flows coupled with a diffusion equation.
We also discuss the variations of the limit equations as the scattering intensity
changes.
Furthermore, when the scattering effect is sufficiently weak, we can obtain the singular limits of the NSF-P1 model with fully general initial data.

\end{abstract}


 \maketitle

\section{ Introduction}

\subsection{The model}
Radiation hydrodynamics is a branch of hydrodynamics in which the moving fluid absorbs and emits electromagnetic radiation.
In radiation hydrodynamics, the absorption or emission of radiation are sufficient to change the pressure of the material, and therefore change its motion; alternatively, the net momentum exchange between radiation and matter  may alter the motion of the matter directly.
The interested reader can refer to \cite{C,MM,P} for more details.

Radiation hydrodynamics mainly concerns with two contents: the propagation of radiation through a fluid and the effect of radiation on fluid flow. Subsequently, based on the governing laws of fluids, the general equations of radiation hydrodynamics  can be written in the following form (see, e.g., \cite{ MM,P})
\begin{equation}\label{RNSF}
\left\{
\begin{aligned}
  & \partial _t \rho  + \nabla \cdot (\rho u) = 0, \\
  & \partial _t \bigg(\rho u+\frac{1}{c^2}F_r\bigg)+\nabla \cdot (\rho u \otimes u+P\mathbb{I}_3+P_r)=\nabla \cdot \Psi(u),\\
  & \partial_t (\rho E+E_r) + \nabla \cdot [(\rho E+P) u+F_r]
     =\nabla \cdot[\Psi(u) u]+\nabla \cdot(\kappa\nabla \Theta).
\end{aligned}
\right.
\end{equation}
Here, $\rho$, $u=(u_1,u_2,u_3)$ and $\Theta$ denote the density, the velocity and the temperature of fluid, respectively.
The viscous stress tensor $\Psi(u)$ is given by
\begin{equation*}
  \Psi(u)=2\mu \mathbb{D}(u)+\lambda (\nabla \cdot u)\mathbb{I}_3,\;\;\;
  \mathbb{D}(u)=\frac{\nabla u+\nabla u ^\top}{2},
\end{equation*}
where $\mu $ and $\lambda$ are viscosity coefficients satisfying $\mu > 0$ and $2\mu + 3\lambda \geq 0$, and $\mathbb{I}_3$ is the $3\times 3$ identity matrix.
The total energy $E$ is given by $E=e+\frac{|u|^2}{2}$, and $e$ denotes the internal energy.
The pressure $P$ and the internal energy $e$ satisfy the perfect gas relations
\begin{equation}\label{gas}
  P=R\rho \Theta \;\;\mathrm{and}\;\; e=c_V\Theta,
\end{equation}
where the constants $R$ and $c_V$ are the generic gas constant and the specific heat at constant volume, respectively.
$\kappa>0$ is the heat conductivity coefficient.
And $E_r$, $F_r$ and $P_r$ denote the radiation energy density, the radiation flux and the radiation pressure tensor, respectively, which are defined by
\begin{equation*}
  E_r=\frac{1}{c}\int_0^\infty \mathrm{d}\nu\int_{\mathbb{S}^2}
  \mathcal{I}(t,x,\nu,\omega)\mathrm{d}\omega,
\end{equation*}
\begin{equation*}
  F_r=\int_0^\infty \mathrm{d}\nu\int_{\mathbb{S}^2}
   \omega \mathcal{I} (t,x,\nu,\omega)\mathrm{d}\omega,
\end{equation*}
and
\begin{equation*}
  P_r=\frac{1}{c}\int_0^\infty \mathrm{d}\nu\int_{\mathbb{S}^2} \omega \otimes \omega \mathcal{I}(t,x,\nu,\omega)\mathrm{d}\omega,
\end{equation*}
where $\mathcal{I}=\mathcal{I}(t,x,\nu,\omega)$ denotes the radiation intensity depending on the frequency $\nu\in(0,\infty)$ and the direction vector $\omega\in \mathbb{S}^2$, and $c>0$ is the  light speed.

To close the radiation hydrodynamics model \eqref{RNSF}, we need to state the governing equation of $\mathcal{I}$ which is a  linear Boltzmann-type equation and takes the form:
\begin{equation}\label{trans-I}
  \frac{1}{c}\partial_t \mathcal{I}+\omega\cdot \nabla \mathcal{I}=S,
\end{equation}
where $S$ stands for the radiative source term.
In this paper, we assume that the radiation fluid is in a state of local thermodynamic equilibrium (see \cite{P}), and then $S$ is defined by
\begin{equation*}
  S=\sigma_a[B(\nu,\Theta)-\mathcal{I}(t,x,\nu,\omega)]+\sigma_s\left(\frac{1}{4\pi}
  \int_{\mathbb{S}^2}\mathcal{I}(t,x,\nu,\omega')\mathrm{d}\omega'
    -\mathcal{I}(t,x,\nu,\omega)\right).
\end{equation*}
The first term in the right-hand of the above equality is the emission-absorption contribution and the second term is the scattering contribution.
$\sigma_a=\sigma_a (\nu,\Theta)\geq 0$ and $\sigma_s=\sigma_s (\nu,\Theta)\geq 0$ denote the absorption coefficient and the scattering coefficient, respectively.
The emission term $B(\nu,\Theta)$ can be taken as the Planck function
\begin{equation*}
  B(\nu,\Theta)=\frac{2h\nu^3}{ c^{2}}\left(e^{h\nu/ k_B\Theta}-1\right)^{-1},
\end{equation*}
where $h$ and $k_B$ are the Planck and Boltzmann constants, respectively.

Putting \eqref{RNSF}-\eqref{trans-I} together, we get the following Navier-Stokes-Fourier-Radiation  (NSF-R) model
\begin{equation}\label{re-RNSF}
\left\{
\begin{aligned}
  & \partial _t \rho  +  u\cdot \nabla \rho+ \rho\nabla \cdot  u = 0, \\
  & \rho (\partial _t u + u\cdot \nabla u) + \nabla P
      =\nabla\cdot \Psi(u) -S_F,\\
  & \rho (\partial _t e + u\cdot \nabla e)+P\nabla \cdot u=
     \nabla\cdot(\kappa \nabla \Theta) +\Psi(u):\nabla u-S_E+u\cdot S_F,  \\
  & \frac{1}{c}\partial_t \mathcal{I}+\omega \cdot \nabla \mathcal{I}=S,
\end{aligned}
\right.
\end{equation}
where
\begin{equation*}
  \nabla\cdot \Psi(u)=\mu \Delta u+(\mu+\lambda)\nabla \nabla \cdot u,
\end{equation*}
\begin{equation*}
   \Psi(u):\nabla u=2\mu |\mathbb{D}(u)|^2+\lambda (\nabla \cdot u)^2,
\end{equation*}
\begin{equation*}
  S_F=\frac{1}{c}\int_0^\infty \mathrm{d}\nu\int_{\mathbb{S}^2}\omega S\mathrm{d}\omega,
\end{equation*}
and
\begin{equation*}
  S_E=\int_0^\infty \mathrm{d}\nu\int_{\mathbb{S}^2}S\mathrm{d}\omega.
\end{equation*}
Here $S_F$ (or $S_E$) characterizes the momentum (or the energy) exchange between the radiation and the matter.

Since the system \eqref{re-RNSF} is very complicated, many  simplified models of it are introduced.
The most  widely-used one is the so-called  Navier-Stokes-Fourier-P1 (NSF-P1) approximation model as described below.
First, if the radiation field is almost isotropic, we can make the assumption   P1 approximation \cite{P}: the radiation intensity $\mathcal{I}$ is given by the first two terms in a spherical harmonic expansion, i.e.,
\begin{equation}\label{ansatz}
  \mathcal{I}(t,x,\nu,\omega)=\frac{1}{4\pi}\mathcal{I}_0(t,x,\nu)+\frac{3}{4\pi}\omega\cdot \mathcal{I}_1(t,x,\nu),
\end{equation}
where the dominant term $\mathcal{I}_0$ and the correction term $\mathcal{I}_1$ are independent of $\omega$.
Plugging \eqref{ansatz} in \eqref{re-RNSF}$_4$, and computing the zero and first order moments with respect to $\omega$ by using the following formulas of solid angle (see \cite{Pe})
\begin{equation*}
  \int_{\mathbb{S}^2}\mathrm{d}\omega=4\pi,\;\;
  \int_{\mathbb{S}^2}\omega\,\mathrm{d}\omega=0,\;\;
  \int_{\mathbb{S}^2}\omega(\omega\cdot A)\mathrm{d}\omega=\frac{4\pi}{3}A,
\end{equation*}
we arrive at
\begin{equation}\label{P1}
\left\{
\begin{aligned}
  & \partial _t \rho  +  u\cdot \nabla \rho+ \rho\nabla \cdot  u = 0, \\
  & \rho (\partial _t u + u\cdot \nabla u) + \nabla P
      =\nabla\cdot \Psi(u)
      +\frac{1}{c}\int_0^\infty (\sigma_a+\sigma_s)\mathcal{I}_1 \mathrm{d}\nu,\\
  & \rho (\partial _t e + u\cdot \nabla e)+P\nabla \cdot u=
     \nabla\cdot(\kappa \nabla \Theta) +\Psi(u):\nabla u\\
  &\qquad
     -\int_0^\infty \sigma_a[4\pi B(\nu, \Theta)-\mathcal{I}_0] \mathrm{d}\nu
     -u\cdot \frac{1}{c}\int_0^\infty (\sigma_a+\sigma_s)\mathcal{I}_1 \mathrm{d}\nu,  \\
  & \frac{1}{c}\partial_t \mathcal{I}_0+\nabla \cdot \mathcal{I}_1
     =\sigma_a[4\pi B(\nu,\Theta)-\mathcal{I}_0],\\
   & \frac{1}{c}\partial_t \mathcal{I}_1+\frac{\nabla \mathcal{I}_0}{3}
     =-(\sigma_a+\sigma_s)\mathcal{I}_1.
\end{aligned}
\right.
\end{equation}

Next, we further assume that the $\sigma_a$ and $\sigma_s$ are independent of the frequency, which is called the gray hypothesis.
Denote $\int^{\infty}_0 \mathcal{I}_0\mathrm{d} \nu$ and $\int^{\infty}_0 \mathcal{I}_1\mathrm{d} \nu$ by $I_0$ and $I_1$, respectively.
Integrating the equations \eqref{P1}$_4$ and \eqref{P1}$_5$ with respect to $\nu$, and noticing that the well-known integration
formula
\begin{equation*}
 \int_0^\infty 4\pi B(\nu,\Theta) \mathrm{d}\nu
  =4\pi \int_0^\infty \frac{2h\nu^3}{ c^{2}}\left(e^{h\nu/ k_B\Theta}-1\right)^{-1} \mathrm{d}\nu=ca_r\Theta^4,
\end{equation*}
where $a_r=\frac{8\pi^5k_B^4}{15h^3c^3}$ is the radiation constant (see \cite{BD,P}),
we obtain the well-known NSF-P1 approximation radiation model:
\begin{equation}\label{NSFP1}
\left\{
\begin{aligned}
  & \partial _t \rho  +  u\cdot \nabla \rho+ \rho\nabla \cdot  u = 0, \\
  & \rho (\partial _t u + u\cdot \nabla u) + \nabla P
      =\nabla\cdot \Psi(u) +\frac{1}{c} (\sigma_a+\sigma_s)I_1,\\
  & \rho (\partial _t e + u\cdot \nabla e)+P\nabla \cdot u=
     \nabla\cdot(\kappa \nabla \Theta)+\Psi(u):\nabla u
     -\sigma_a(ca_r \Theta^4-I_0)-\frac{1}{c} (\sigma_a+\sigma_s)u\cdot I_1,  \\
  & \frac{1}{c}\partial_t I_0+\nabla \cdot I_1
     =\sigma_a[ca_r \Theta^4-I_0],\\
   & \frac{1}{c}\partial_t I_1+\frac{\nabla I_0}{3}
     =-(\sigma_a+\sigma_s)I_1.
\end{aligned}
\right.
\end{equation}

In order to identify the relevant singular limit regime, we reformulate \eqref{NSFP1} into the following dimensionless form, which only retains the useful parameters and ignores the influence of other parameters (we give the detailed derivation
in the  Appendix),
\begin{equation}\label{NSFP1-Ma}
\left\{
\begin{aligned}
  & \partial _t \rho  +  u\cdot \nabla \rho+ \rho\nabla \cdot  u = 0, \\
  & \rho \partial _t u + \rho u\cdot \nabla u + \frac{1}{ \mathrm{Ma^2}}\nabla P
      =\nabla\cdot \Psi(u)
        +\frac{\mathcal{P}\mathcal{L}}{\mathrm{Ma^2}}
        (1+\mathcal{L}_s)I_1,\\
  & \rho \partial _t \Theta + \rho u\cdot \nabla \Theta+ P\nabla \cdot u =\kappa\Delta \Theta
      +\mathrm{Ma}^2\Psi(u):\nabla u
      +\mathcal{PCL}(I_0-\Theta^4)-\mathcal{PL}(1+\mathcal{L}_s)I_1\cdot u,  \\
  & \frac{1}{\mathcal{C}} \partial_t I_0
    +\nabla \cdot  I_1=\mathcal{L}(\Theta^4-I_0),\\
  & \frac{1}{\mathcal{C}} \partial_t I_1
     +\frac{1}{3}\nabla I_0=-\mathcal{L}(1+\mathcal{L}_s)I_1.
\end{aligned}
\right.
\end{equation}
Here, $\mathrm{M}_a$ and $\mathcal{C}$ are Mach and ``infrarelativistic" numbers, and $\mathcal{P}$, $\mathcal{L}$ and $\mathcal{L}_s$ are various dimensionless numbers corresponding to the radiation, see the Appendix for the expressions and physical meanings of these dimensionless numbers.
We would like to point out that $\mathcal{L}$ and $\mathcal{L}\mathcal{L}_s$ measure the absorption strength and the scattering strength, respectively.

\subsection{Previous results on singular limits of the models of radiation hydrodynamics}
Singular limits of the models of radiation hydrodynamics mainly involve two types: diffusion limit and low Mach number limit. These two regimes are not determined by a single parameter, but by multiple parameters.

The diffusion limits of radiation fluid models have been introduced by Pomraning \cite{P}.
From a physical point of view, when the mean free-path of a photon is small enough, the radiative transfer equation can be approximated by a diffusion equation.
This asymptotic behavior has been studied formally and numerically by Lowrie, Morel and  Hittinger \cite{LMH} and Buet and Despres \cite{BD}.
The diffusion limits of the models of radiation hydrodynamics involves two cases.
To be specific, when the emission-absorption effect is dominant, the corresponding case is called equilibrium diffusion regime
(i.e. $\mathcal{L}=O(\e^{-1})$, $\mathcal{L}_s=O(\e^2)$, $\mathrm{Ma}=O(1)$ and  $\mathcal{C}=O(\e^{-1})$).
On the contrary, when the scattering effect is dominant, the corresponding case is called non-equilibrium diffusion regime
(i.e. $\mathcal{L}=O(\e)$, $\mathcal{L}_s=O(\e^{-2})$, $\mathrm{Ma}=O(1)$ and $\mathcal{C}=O(\e^{-1})$).
The rigorous proof of both equilibrium and non-equilibrium diffusion limits for the NSF-R model \eqref{re-RNSF} in the framework of weak solutions has been given by Ducomet and  Ne\v{c}asov\'{a} \cite{DN2,DN3} thorough using the relative entropy method
for well-prepared initial data.
Later, Danchin and Ducomet \cite{DD1} established the existence of global strong solutions to the isentropic Navier-Stokes-P1 model with small enough initial data in critical regularity spaces, and the influence of absorption and scattering coefficients on the equilibrium and non-equilibrium diffusion limits was discussed in more detail.

In many combustion phenomena, the characteristic speed of flow is very small, while the characteristic temperature is very large and the effect of thermal radiation cannot be ignored.
In this situation, it is meaningful to consider the low Mach number limit of the models of radiation hydrodynamics (i.e. $\mathcal{L}=\mathcal{L}_s=O(1)$, $\mathrm{Ma}=O(\e)$ and $\mathcal{C}=O(\e^{-1})$).
In the framework of weak solutions for well-prepared  initial data, Ducomet and Ne\v{c}asov\'{a} \cite{DN1} investigated the low Mach number limit for the NSF-R model \eqref{re-RNSF} by using the relative entropy method, and proved the convergence toward the incompressible Navier-Stokes system coupled to a system of two stationary transport equations.
In the framework of classical solutions, Danchin and Ducomet \cite{DD} studied the low Mach number limit of isentropic Navier-Stokes-P1 model with well-prepared initial data, and showed the convergence to the incompressible Navier-Stokes equations.
Fan, Li and Nakamura \cite{FLN} proved that the NSF-P1 model \eqref{NSFP1-Ma} also converges to the incompressible Navier-Stokes equations for well-prepared initial data and small variation of temperature.
In addition, we refer to \cite{DST,TS,TSGKS} for the physical meanings and numerical simulation of radiation models at low Mach number.

We note that, in all above results, the temperature is restricted to a small variation or remains constant. At the same time, the initial data of the corresponding model are well-prepared.
However, in engineering applications and mathematical theory, it is more important to consider large temperature variations and general initial data (also known as ill-prepared initial data).
Recently, Jiang, Ju and Liao \cite{JJL} investigated the combined non-equilibrium diffusion and low Mach number limits of the Euler-P1 approximation model (i.e. $\mathcal{L}=O(\e)$, $\mathcal{L}_s=O(\e^{-2})$, $\mathrm{Ma}=O(\e)$ and $\mathcal{C}=O(\e^{-1})$) with large temperature variations,
but they also need the well-prepared initial data and ignore heat conduction.

The main purpose of this paper is to study the combined non-equilibrium diffusion and low Mach number limits of the NSF-P1 model \eqref{NSFP1-Ma} with \emph{heat conduction} and \emph{large temperature variations}.
In particular, we broaden the range of initial data from well-prepared initial data to \emph{partially general ones} (see Remark \ref{haochuzhi} below) and hence improve the results in \cite{FLN,JJL}.
Furthermore, we introduce a parameter $\delta \in [0,2]$ to describe the scattering intensity ($\mathcal{L}_s=\e^{-\delta}$) and discuss the variations of the limit equations as the scattering intensity changes. We find that, with the weakening of scattering intensity  ($\delta$ from 2 to 0), the ``diffusion property" of the dominant term  $I_0$ gradually weakens and the ``importance" of the correction term $I_1$ is gradually increasing.
When the scattering effect is sufficiently weak  ($\delta=0$), we establish the limit system under \emph{fully general initial data}.

\subsection{Our results}
Taking
\begin{equation*}
\mathrm{Ma}=\mathcal{P}=\mathcal{L}=\e,\;\mathcal{C}=\e^{-1},\;\mathrm{and}\;\;
\mathcal{L}_s=\e^{-\delta},\;\delta\in [0,2]
\end{equation*}
in \eqref{NSFP1-Ma} yields
\begin{equation}\label{guodu}
\left\{
\begin{aligned}
  & \partial _t \rho  +  u\cdot \nabla \rho+ \rho\nabla \cdot  u = 0, \\
  & \rho( \partial _t u + u \cdot \nabla u) + \frac{1}{\e^2}\nabla P
      =\nabla\cdot \Psi(u)+\bigg(1+\frac{1}{\e^\delta}\bigg)I_1,\\
  & \rho( \partial_t \Theta +  u\cdot \nabla \Theta)+ P\nabla \cdot u=
       \kappa\Delta \Theta +\epsilon^2 \Psi(u):\nabla u
      +\epsilon(I_0-\Theta^4)-(\epsilon^2 +\e^{2-\delta}) I_1\cdot u,  \\
  & \epsilon \partial_t I_0
    +\nabla \cdot  I_1=\epsilon(\Theta^4-I_0),\\
  & \epsilon \partial_t I_1
     +\frac{1}{3}\nabla I_0=-\epsilon\bigg(1+\frac{1}{\e^{\delta}}\bigg) I_1.
\end{aligned}
\right.
\end{equation}
We introduce the scalings of pressure and temperature as
\begin{equation}\label{ptheta}
  P=e^{\epsilon p^\epsilon} \;\; {\rm and}\;\; \Theta=e^{ \theta^\epsilon},
\end{equation}
which mean that $P \sim 1+\e p^\e$ and $\Theta \sim 1+\theta^\e$, and imply  the large variation of temperature.
Putting these scalings into \eqref{guodu} and using the dimensionless relation $P=\rho\Theta$, we rewrite \eqref{guodu} as
\begin{equation}\label{nondim}
\left\{
\begin{aligned}
  &  \partial_t p^\epsilon+u^\epsilon\cdot \nabla p^\epsilon
     +\frac{2}{\epsilon}\nabla \cdot u^\e = e^{-\e p^\epsilon}
     \bigg[\frac{\kappa}{\epsilon} \Delta e^{\theta^\e}
      +\e \Psi(u^\e):\nabla u^\e + (I^\e_0-e^{4\theta^\e})
      - \bigg(\e+\frac{1}{\e^{\delta-1}}\bigg) I^\e_1\cdot u^\e\bigg], \\
  &  e^{-\theta^\e} (\partial_t u^\epsilon + u^\epsilon\cdot \nabla u^\epsilon )
       + \frac{\nabla p^\epsilon}{\epsilon}
     =e^{-\e p^\e}\bigg[\nabla\cdot \Psi(u^\e)
        +\bigg(1+\frac{1}{\e^\delta}\bigg)I^\e_1\bigg],\\
  & \partial_t \theta^\e + u^\e\cdot \nabla \theta^\e+ \nabla \cdot u^\e
    = e^{-\e p^\epsilon}\big[ \kappa \Delta e^{\theta^\e}
    +\epsilon^2\Psi(u^\e):\nabla u^\e
      +\epsilon(I^\e_0-e^{4\theta^\e})- (\epsilon^2+\e^{2-\delta}) I^\epsilon_1\cdot u^\epsilon\big],  \\
  &  \partial_t I^\epsilon_0
    +\frac{\nabla \cdot I^\epsilon_1}{\epsilon}=e^{4\theta^\e}-I^\epsilon_0,\\
  & \partial_t I^\epsilon_1
     +\frac{\nabla I^\epsilon_0}{3\epsilon}=-\bigg(1+\frac{1}{\e^\delta}\bigg)I^\epsilon_1,
\end{aligned}
\right.
\end{equation}
where we have added the superscript ``$\epsilon $" on the unknowns $(p^\e,u^\e,\theta^{\e},I^\e_0,I^\e_1)$ to emphasize the dependence on $\epsilon$.
We supply the system \eqref{nondim} with the initial data
\begin{equation}\label{initial data}
  (p^\e,u^\e,\theta^{\e},I^\e_0,I^\e_1)|_{t=0}
  =(p^\e_0(x),u^\e_0(x),\theta^{\e}_0(x),I^\e_{00}(x),I^\e_{10}(x)),\;\;\;x\in \mathbb{R}^3.
\end{equation}

Define the weighted norms:
\begin{equation*}
  \|v\|_{s,\epsilon,T}:= \sup\limits_{0\leq t\leq T}\|v(t)\|_{s,\epsilon},\;\;\;
    \|v(t)\|_{s,\epsilon}:= \sum\limits_{k=0}^s\|(\epsilon\partial_t)^k v(t)\|_{s-k},
\end{equation*}
and
\begin{equation*}
  \interleave v\interleave_{s,\epsilon,T}
  := \sup\limits_{0\leq t\leq T} \interleave v(t) \interleave_{s,\epsilon},\;\;\;
     \interleave v(t) \interleave_{s,\epsilon}:= \sum\limits_{k=0}^s\|\epsilon^{[k-1]^+}\partial_t^k v(t)\|_{s-k},
\end{equation*}
where $\|\cdot\|_s $ denote the norm of $H^s(\mathbb{R}^3)$ and $ [k-1]^+=\max\{k-1,0\}$.

We first consider the case $\delta\in(0,2]$, and state the uniform existence of the local solutions as follows.
\begin{theorem}\label{uniform existence theorem}
Let $\delta\in(0,2]$ and $s\geq 4$ be an integer. Assume that the initial data \eqref{initial data} satisfy
\begin{equation}\label{initial condition}
  \|(p^\e_0,u^\e_0) \|_s
  +\|(\e p^\e_0,\e u^\e_0,\theta^{\e}_0-\theta_c) \|_{s+1}
  + \interleave (I^\e_{00}-I_c, I^\e_{10}) \interleave_{s+1,\e}
  \leq M_0,
\end{equation}
for the three positive constants $\theta_c$, $I_c$ and $M_0$ independent of $\epsilon$, and $\theta_c$ and $I_c$  satisfying $I_c=e^{4\theta_c}$.
Then there exist constants $\epsilon_0\in (0,1]$ and $T_0>0$ such that, for all $\epsilon \leq \epsilon_0$, the Cauchy problem \eqref{nondim} and \eqref{initial data} has a unique smooth solution $(p^\e,u^\e,\theta^{\e},I^\e_0,I^\e_1)$
on $[0,T_0]$ satisfying
\begin{align}
  &\|(p^\e,u^\e)\|_{s,\e,T_{0}}
     +\|(\e p^\e,\e u^\e, \theta^\e-\theta_c)\|_{s+1,\e,T_0}
     + \interleave (I^\e_0-I_c,I^\e_1) \interleave_{s+1,\e,T_0}  \notag\\
    & \qquad+ \bigg(\int_0^{T_0} \|\nabla(p^\e,u^\e)(t)\|^2_{s,\e}
      +\|\nabla(\e u^\e,\theta^\e-\theta_c)(t)\|^2_{s+1,\e}
      \mathrm{d}t\bigg)^{\frac{1}{2}}  \notag\\
    &\qquad\qquad+ \bigg(\int_0^{T_0}
      \bigg(1+\frac{1}{\e^\delta}\bigg)
      \interleave\!(I^\e_0-I_c,I^\e_1)(t)\! \interleave^{2}_{s+1,\epsilon}
      \mathrm{d}t\bigg)^{\frac{1}{2}} \leq M_1, \label{uniform estimate}
\end{align}
where the constant $M_1>0$ depends only on $\theta_c$, $I_c$, $M_0$ and $T_0$.
\end{theorem}
\begin{remark}
   In the assumption \eqref{initial condition}, $ \partial_t I^\e_{00}$ is indeed defined by $
    \partial_t I^\epsilon_{00}=
    -\frac{\nabla \cdot I^\epsilon_{10}}{\epsilon}+e^{4\theta^\e_0}-I^\epsilon_{00}$ through the  equation \eqref{nondim}$_4$.
    And $ \partial_t I^\e_{10}$  is defined by an analogous way.
\end{remark}

\begin{remark}\label{haochuzhi}
In singular limits problems, well-prepared initial data means that there is no initial layer.
On the contrary, general initial data always lead to the generation of initial layer.
We call \eqref{initial condition} the partially general initial data condition since
only the boundedness of $(\partial_t I^\e_{00},\partial_t I^\e_{10})$ is required, and
there are no additional assumptions on
$(\partial_t p^\e_0,\partial_t u^\e_0,\partial_t \theta^\e_0)$.
Compared with the initial data conditions stated in \cite{JJL}, which include the boundedness of  $(\partial_t p_0^\e,\partial_t u_0^\e, \partial_t \theta_0^\e,\partial_t I_{00}^\e, \partial_t I^\e_{10})$  and of $(\partial^2_t I_{00}^\e, \partial^2_t I^\e_{10})$, we relax the restriction on the initial data to a great extent.
\end{remark}


The convergence results for the parameter $\delta\in (0,2]$ read as follows.
\begin{theorem}\label{convergence thm}
Suppose that the assumptions in Theorem \ref{uniform existence theorem} hold. Assume further that the initial data \eqref{initial data} satisfy
\begin{equation*}
  ( p^\e_0, u^\e_0,\theta^{\e}_0-\theta_c,I^\e_{00}-I_c,I^\e_{10})
  \rightarrow
  \bigg(-\frac{\bar{I}_{00}-I_c}{3}, \bar{u}_0,\bar{\theta}_0-\theta_c,\bar{I}_{00}-I_c,\bar{I}_{10}\bigg)\;\;
  \mathrm{in}\;\;H^s(\mathbb{R}^3)
\end{equation*}
as $\e \rightarrow 0$, and $\theta^\e_0$ decays at infinity in the sense that
\begin{equation}\label{decay theta}
  |\theta^\e_0(x)-\theta_c|\leq c_0|x|^{-1-\sigma},\quad
  |\nabla\theta^\e_0(x)|\leq c_0|x|^{-2-\sigma},
\end{equation}
where $c_0$ and $\sigma$ are given positive constants.
Then the solution of the Cauchy problem \eqref{nondim} and \eqref{initial data} satisfies
\begin{equation*}
  ( p^\e, u^\e,\theta^{\e}-\theta_c,I^\e_{0}-I_c,I^\e_{1})\rightarrow
  \bigg(-\frac{\bar{I}_{0}-I_c}{3},\bar{u},\bar{\theta}-\theta_c,\bar{I}_0-I_c,\bar{I}_1\bigg)
\end{equation*}
weakly$-*$ in $L^{\infty}(0,T_0;H^s(\mathbb{R}^3))$ and strongly in $L^2 (0,T_0;H^{s'}_{\mathrm{loc}}(\mathbb{R}^3))$ for any $s'\in[0,s)$.
Moreover,
\begin{itemize}
\item when $\delta = 2$, then $\bar{I}_{1}=0$, and there exists some function $\pi_1 \in C([0,T_0];H^{s}(\mathbb{R}^3))$ such that $(\bar{u},\bar{\theta}, \bar{I}_0)$ satisfies the system of low Mach number heat-conducting viscous flows coupled with a diffusion equation
  \begin{equation}\label{limit eq}
  \left\{
  \begin{aligned}
  &   2\nabla \cdot\bar{u}=\kappa \Delta e^{\bar{\theta}},\\
  &  e^{-\bar{\theta}} (\partial_t \bar{u} + \bar{u}\cdot \nabla \bar{u} )
       +\nabla \pi_1=\nabla\cdot \Psi(\bar{u}),\\
  & \partial_t \bar{\theta} + \bar{u}\cdot \nabla \bar{\theta}+ \nabla \cdot \bar{u}
    =   \kappa \Delta e^{\bar{\theta}},\\
  & \partial_t \bar{I}_0-\frac{1}{3}\Delta \bar{I}_0+\bar{I}_0= e^{4\bar{\theta}},
  \end{aligned}
  \right.
  \end{equation}
  with the initial data
  $(\bar{u},\bar{\theta},\bar{I}_{0})|_{t=0}=(\bar{w}_0,\bar{\theta}_0,\bar{I}_{00})$,
  where $\bar{w}_0$ is determined by
  \begin{equation}\label{in2}
  2\nabla \cdot\bar{w}_0=\kappa \Delta e^{\bar{\theta}_0},\;\;
  \nabla\times (e^{-\bar{\theta}_0}\bar{w}_0)=\nabla\times (e^{-\bar{\theta}_0}\bar{u}_0);
  \end{equation}
\item  when $\delta\in(1,2)$, then $\bar{I}_{1}=0$, $(\bar{u},\bar{\theta})$ solves the equations \eqref{limit eq}$_1$-\eqref{limit eq}$_3$ with the initial data  $(\bar{u},\bar{\theta})|_{t=0}=(\bar{w}_0,\bar{\theta}_0)$ satisfying \eqref{in2}, and $\bar{I}_0$  satisfies a Laplace equation
  \begin{equation*}
   \Delta \bar{I}_0=0;
  \end{equation*}
\item  when $\delta=1$, then $(\bar{u},\bar{\theta})$ solves the equations
   \eqref{limit eq}$_1$-\eqref{limit eq}$_3$ with the initial data  $(\bar{u},\bar{\theta})|_{t=0}=(\bar{w}_0,\bar{\theta}_0)$ satisfying \eqref{in2},
   and $(\bar{I}_0,\bar{I}_1)$  satisfies
  \begin{equation*}
    \nabla \bar{I}_0=-3\bar{I}_1,\;\;\; \nabla \cdot \bar{I}_1=0;
  \end{equation*}
\item  when $\delta\in(0,1)$, then $(\bar{u},\bar{\theta})$ solves the equations
   \eqref{limit eq}$_1$-\eqref{limit eq}$_3$ with the initial data  $(\bar{u},\bar{\theta})|_{t=0}=(\bar{w}_0,\bar{\theta}_0)$ satisfying \eqref{in2},
   and $(\bar{I}_0,\bar{I}_1)$ satisfies
  \begin{equation*}
   \nabla \bar{I}_0=0,\;\;\; \nabla \cdot \bar{I}_1=0.
  \end{equation*}
  \end{itemize}
\end{theorem}

\begin{remark}\label{re-1.2}
The difference between the above four limit equations is reflected in the characterization of $(\bar{I}_0,\bar{I}_1)$.
Generally speaking, the convergence result of the first case is also called the non-equilibrium diffusion limit at low Mach number in radiation hydrodynamics, see \cite{JJL}.
With the weakening of scattering intensity, the ``diffusion property" of the dominant term $I_0$ gradually weakens and the ``importance" of the correction term $I_1$ is gradually increasing.
\end{remark}

\begin{remark}
Recall the scalings \eqref{ptheta}.
Taking $\e \rightarrow 0$ yields
\begin{equation*}
\bar{\Theta}= e^{\bar{\theta}}  \;\; \mathrm{and}\;\;\bar{\Theta}\bar{\rho}=1.
\end{equation*}
Then we can reformulate \eqref{limit eq}$_1$-\eqref{limit eq}$_3$ into the following low Mach number inhomogeneous Navier-Stokes equations
\begin{equation}\label{re-rho}
  \left\{
  \begin{aligned}
  &  \partial_t \bar{\rho}+\bar{u}\cdot\nabla \bar{\rho}+\bar{\rho}\nabla\cdot\bar{u}=0,\\
  &  \bar{\rho}(\partial_t \bar{u} + \bar{u}\cdot \nabla \bar{u} )
       +\nabla \pi_1=\nabla\cdot \Psi(\bar{u}),\\
  &   2\nabla \cdot\bar{u}
      =\nabla \cdot \bigg(\kappa\nabla\bigg(\frac{1}{\bar{\rho}}\bigg)\bigg).\\
  \end{aligned}
  \right.
\end{equation}
The formal derivation from the compressible Navier-Stokes-Fourier equations to the equations \eqref{re-rho}, as Mach number tends to zero, is given in P.-L. Lions' famous book \cite{L}.
Especially, if we take $\kappa=0$ in \eqref{re-rho}, we obtain the so-called inhomogeneous incompressible Navier-Stokes equations.
\end{remark}

Now, let's give some comments on the proofs of Theorems
\ref{uniform existence theorem}-\ref{convergence thm}.
The uniform estimates of the solution $(p^\e,u^\e,\theta^{\e},I^\e_0,I^\e_1)$ are the main part of Theorem \ref{uniform existence theorem}.
There are two difficulties in getting uniform estimates.
The first one is caused by the radiation pressure and the effect of heat conduction, which lead to more complex singular terms in the equation \eqref{nondim} destroying the symmetric singular structure of this system.
Thus, the classical theory developed by Klainerman and Majda \cite{KM} is not applicable.
To surrounding this difficulty, we construct auxiliary equations of
$(\e p^\e-(\theta^\e-\theta_c),\e u^\e,\theta^\e-\theta_c,I^\e_0-I_c,I^\e_1)$ (see \eqref{nondim-hat} below), which own a symmetric singular structure structure.
At the same time, we can use these auxiliary equations to get higher order derivative estimates of $(\e p^\e,\e u^\e,\theta^\e-\theta_c)$, which play a key role in the whole proof.

When we establish the estimates of $(p^\e,u^\e)$ by energy method, we encounter the second difficulty that the spatial-temporal mixed derivative estimates of $(p^\e,u^\e)$ can not match with the order of $\e$.
This difficulty is caused by the large temperature variation, which leads to the generation of unbounded terms during energy estimates (The same problem also appears in the non-isentropic Euler equations and in the ideal non-isentropic magnetohydrodynamic equations, see \cite{MS,LZ} for more explanations).
Our strategy here is to divide $\|( p^\e, u^\e)\|_{s,\e}$ into three parts:
\begin{equation*}
  \sum\limits_{k=0}^{s}\|(\e\partial_t)^k (p^\e,u^\e)\|,\;\;
  \|(\nabla p^\e,\nabla \cdot u^\e)\|_{s-1,\e}\;\;
   \mathrm{and}\;\; \|\nabla \times u^\e\|_{s-1,\e}.
\end{equation*}
To obtain the estimate of the first part, we introduce the equivalent pressure $\tilde{p}$ and velocity $\tilde{u}$ as
\begin{equation*}
  \tilde{p}=p^\e+\frac{e^{-\e p^\e}(I^\e_0-I_c)}{3}\;\;\mathrm{and}\;\;
  \tilde{u}=2u^\e-\kappa e^{-\epsilon p^\e+\theta^\e}\nabla\theta^\e,
\end{equation*}
and construct the equations of $(\tilde{p},\tilde{u})$ (see \eqref{pt-ut} below).
We point out that the variable $\tilde{p}$ consists of two parts: the fluid pressure and the radiation pressure, and the idea of introducing $\tilde{p}$ is based on the mechanical effect of radiation.
Then we obtain the estimate of $\sum\limits_{k=0}^{s}\|(\e\partial_t)^k (\tilde{p},\tilde{u})\|$, which is equivalent to the first part.
The boundedness of the second part is indirectly obtained by using the structure of the equations \eqref{nondim} and the above obtained estimates.
Finally, we use energy method to get the estimate of $\|\nabla \times(e^{-\theta^\e} u^\e)\|_{s-1,\e}$ and then the estimate of the third part follows immediately.

Once the uniform existence of the solutions have been established, we are in a position to show the convergence results (Theorem \ref{convergence thm}).
Since we consider here the partial general initial data condition \eqref{initial condition}, which doesn't give us the convergence of $(p^\e, u^\e)$ directly.
Hence, we divide the solution into slow components $(e^{-\theta^\e}u^\e,\theta^\e,I^\e_0,I^\e_1)$ and fast components $(\tilde{p},\tilde{u})$.
The compactness of slow components is obtained from the above uniform estimates and Aubin-Lions Lemma \cite{Simon}.
The method we use to get the convergence of fast components is based on the local energy decay of acoustic wave equations, which is developed by M\'{e}tivier and Schochet \cite{MS} on the non-isentropic Euler equations with general initial data, see also \cite{A,JJLX,LST} for further extensions on the full Navier-Stokes equations and the full magnetohydrodynamic equations.


Below we  consider the case $\delta=0$, which means that both the scattering intensity ($\mathcal{L}\mathcal{L}_s=\e$) and the absorption intensity ($\mathcal{L}=\e$) are sufficiently weak. In this situation,
we can establish the uniform existence of the solutions and the convergence to the corresponding limit equations with fully general initial data.
More precisely, the results read as follows.
\begin{theorem}\label{Ls 0 thm}
Let $\delta=0$ and $s\geq 4$ be an integer.
Assume that the initial data \eqref{initial data} satisfy
\begin{equation}\label{initial condition Ls0}
  \|(p^\e_0,u^\e_0) \|_s
  +\|(\e p^\e_0,\e u^\e_0,\theta^{\e}_0-\theta_c,I^\e_{00}-I_c,I^\e_{10}) \|_{s+1}
  \leq M'_0,
\end{equation}
for a positive constant $M'_0$ independent of $\epsilon$.
Then there exist constants $\epsilon'_0\in (0,1]$ and $T'_0>0$ such that, for all $\epsilon \leq \epsilon'_0$, the Cauchy problem \eqref{nondim} and \eqref{initial data} has a unique smooth solution $(p^\e,u^\e,\theta^{\e},I^\e_0,I^\e_1)$
on $[0,T^\prime_0]$ satisfying
\begin{align}
   &\|(p^\e,u^\e)\|_{s,\e,T'_0}
     +\|(\e p^\e,\e u^\e, \theta^\e-\theta_c,I^\e_0-I_c,I^\e_1)\|_{s+1,\e,T'_0} \notag\\
   &\qquad +  \bigg(\int_0^{T'_0} \|\nabla(p^\e,u^\e)(t)\|^2_{s,\e}
      +\|\nabla(\e u^\e,\theta^\e-\theta_c)(t)\|^2_{s+1,\e}
      +\|(I^\e_0-I_c,I^\e_1)(t)\|^2_{s+1,\epsilon}
      \mathrm{d}t\bigg)^{\frac{1}{2}}
      \leq  M'_1,
      \label{Ls-uniform estimate}
\end{align}
where the constant $M'_1>0$ depends only on $\theta_c$, $I_c$, $M'_0$ and $T'_0$.
Furthermore, if we assume that the initial data \eqref{initial data} satisfy \eqref{decay theta} and
\begin{equation*}
  ( p^\e_0, u^\e_0,\theta^{\e}_0-\theta_c,I^\e_{00}-I_c,I^\e_{10})
  \rightarrow
  (0, \bar{u}_0,\bar{\theta}_0-\theta_c,0,\bar{I}_{10})\;\;
  \mathrm{in}\;\;H^s(\mathbb{R}^3)
\end{equation*}
as $\e \rightarrow 0$, then the solution of Cauchy problem \eqref{nondim} and \eqref{initial data} satisfies
\begin{equation*}
  ( p^\e, u^\e,\theta^{\e}-\theta_c,I^\e_{0}-I_c,I^\e_{1})\rightarrow
  (0,\bar{u},\bar{\theta}-\theta_c,0,\bar{I}_1)
\end{equation*}
weakly$-*$ in $L^{\infty}(0,T_0';H^s(\mathbb{R}^3))$ and strongly in $L^2 (0,T_0';H^{s'}_{\mathrm{loc}}(\mathbb{R}^3))$ for any $s'\in[0,s)$,
where $(\bar{u},\bar{\theta},\bar{I}_1)$ solves the following equations
\begin{equation}\label{limit eq-3}
  \left\{
  \begin{aligned}
  &    2\nabla \cdot\bar{u}=\kappa \Delta e^{\bar{\theta}},\\
  &  e^{-\bar{\theta}} (\partial_t \bar{u} + \bar{u}\cdot \nabla \bar{u} )
       +\nabla \pi_2=\nabla\cdot \Psi(\bar{u}),\\
  & \partial_t \bar{\theta} + \bar{u}\cdot \nabla \bar{\theta}+ \nabla \cdot \bar{u}
    =   \kappa \Delta e^{\bar{\theta}},\\
  & \partial_t \bar{I}_1+\nabla \pi_3= -2\bar{I}_1,\;\; \nabla \cdot \bar{I}_1=0,
  \end{aligned}
  \right.
\end{equation}
for some functions $\pi_2$, $\pi_3\in C([0,T'_0];H^s(\mathbb{R}^3))$,
and the initial data $(\bar{u},\bar{\theta}, \bar{I}_1)|_{t=0}=(\bar{w}_0,\bar{\theta}_0,\bar{I}_{10})$ satisfy \eqref{in2} and $\nabla \cdot \bar{I}_{10}=0$.
\end{theorem}

\begin{remark}
Compared with the initial condition \eqref{initial condition} in
Theorem \ref{uniform existence theorem}, there are no restrictions on $(\partial_t I^\e_{00},\partial_t I^\e_{10})$ in \eqref{initial condition Ls0},
which means that \eqref{initial condition Ls0} is the fully general initial data condition.
\end{remark}

\begin{remark}
From a physical point of view, $I_0$ is dominant and $I_1$ represents the first order anisotropy correction to $I_0$.
The disappearance of $\bar{I}_0$ in the limit system \eqref{limit eq-3} is therefore surprising.
From \cite{P,Pe}, we found the reason is that when the absorption intensity and scattering intensity are not strong enough, the radiation field is not almost isotropic and then the P1 approximation \eqref{ansatz} can not approximate the radiation intensity well.
Nevertheless, from a mathematical point of view, we still believe that our results of the case $\delta=0$ with fully general initial data are very meaningful.
\end{remark}

\begin{remark}
Based on the differences of the limit systems of the five  cases: $\delta=2$, $\delta\in(1,2)$, $\delta=1$, $\delta\in(0,1)$ and $\delta=0$, we find that,
as $\delta$ decreases, the accuracy of the P1 approximation decreases.
\end{remark}

Theorem \ref{Ls 0 thm} contains two parts, the uniform existence and convergence of the classical solutions when $\delta=0$.
Although the initial data \eqref{initial condition Ls0} are general, the method we used to prove Theorem \ref{uniform existence theorem} is still valid for proving the uniform existence part of Theorem \ref{Ls 0 thm}.
Compared with Theorem \ref{convergence thm}, since the compactness of $(I^\e_0,I^\e_1)$ cannot be derived from the general initial data \eqref{initial condition Ls0},
the key to proving the convergence part of Theorem \ref{Ls 0 thm} is to obtain the convergence of $(I^\e_0,I^\e_1)$ in other ways.
We remark that, when $\delta=0$, the singular structure of equations satisfied by $(I^\e_0,I^\e_1)$ is similar to that of the low Mach number regime for compressible fluid, which plays a fundamental role in our analysis.
According to this structure, we first show the convergence of $\nabla \times I^\e_1$.
And then we construct the wave equation satisfied by $(\nabla I^\e_0, \nabla \cdot I^\e_1)$ and acquire the convergence of $(\nabla I^\e_0, \nabla \cdot I^\e_1)$ by using the local energy decay of wave equations.

\subsection*{Notations}

For a multi-index $\alpha=(\alpha_1,\alpha_2,\alpha_3)$, we denote $\partial_x^\alpha=\partial_{x_1}^{\alpha_1}\partial_{x_2}^{\alpha_2}\partial_{x_3}^{\alpha_3}$ and $|\alpha|=\alpha_1+\alpha_2+\alpha_3$.
For an integer $k$ and a multi-index $\alpha$, we denote
$D^{k,\alpha}=(\e\partial_t)^k\partial_x^\alpha$ and
$\dot{D}^{k,\alpha}=\e^{[k-1]^+}\partial_t^k\partial_x^\alpha$.
The symbol $D^{k,i}$ (or $D^i$) denotes the summation of all $D^{k,\alpha}$ with $|\alpha|=i$  (or $k+|\alpha|=i$).

We use $L^2(\mathbb{R}^3)$ to denote the space of square integrable functions on $\mathbb{R}^3$ with the norm $\|\cdot\|$.
The inner product in $L^2(\mathbb{R}^3)$ is denoted by $\langle\cdot,\cdot\rangle$.
$L^{\infty}(\mathbb{R}^3)$ is the space of essentially bounded functions on $\mathbb{R}^3$ with the norm $\| \cdot \|_{L^{\infty}}$.
$H^s (\mathbb{R}^3)$ denotes the standard Sobolev spaces  $W^{s,2}(\mathbb{R}^3)$ with the norm $\| \cdot \|_s$.
Furthermore, we denote by $C^i([0,T];H^s (\mathbb{R}^3))$ the space of $i$-th times continuously differentiable functions on $[0,T]$ taking values in $H^s (\mathbb{R}^3)$.

We use $C_0(\cdot)$ and $C(\cdot)$ to denote two  positive increasing polynomial functions from $[0,\infty)$ to $[0,\infty)$ independent of $\epsilon$, which may vary from line to line.
The notation $A\lesssim B$ means that $A \leq CB$ holds for some positive constant $C$ independent of $\epsilon$.

\vskip 3mm
The rest of this paper is arranged as follows.
In the next section, we give some basic facts and inequalities.
We establish the uniform estimates of the solutions and prove Theorem \ref{uniform existence theorem} in Section \ref{S.3}.
In Section \ref{S.4}, we study the dispersive estimates on acoustic wave equations satisfied by the fast components and prove Theorem \ref{convergence thm}.
In section \ref{S.5}, we give the sketch of proof to Theorem \ref{Ls 0 thm}.
Finally, the dimensional analysis of \eqref{NSFP1} is given in the Appendix.

\section{Preliminaries}\label{S.2}

We first recall the results on Moser-type calculus inequalities and the estimate of composite functions in Sobolev spaces.
\begin{lemma}[\cite{KM,M}]
Let $s\in \mathbb{N}$.
Assume that $u,\,v\in H^s(\mathbb{R}^3)\cap L^{\infty}(\mathbb{R}^3)$.
Then for any $\alpha$ with $1\leq|\alpha|\leq s$, we have
\begin{equation*}
  \|\partial_x^\alpha(uv)\| \lesssim \|u\|_{L^{\infty}}\|\partial_x^s v\|
    +\|v\|_{L^{\infty}}\|\partial_x^s u\|.
\end{equation*}
Assume further that $\nabla u\in L^{\infty}(\mathbb{R}^3)$, then
\begin{equation*}
  \|[\partial_x^\alpha, u]v\|
   \lesssim \|\nabla u\|_{L^{\infty}}\|\partial_x^{s-1} v\|
    +\|v\|_{L^{\infty}}\|\partial_x^s u\|,
\end{equation*}
where $[\partial_x^\alpha,u]v=\partial_x^\alpha (uv)-u \partial_x^\alpha v$.
Moreover, if $s>\frac{5}{2}$, it holds that
\begin{equation*}
\begin{aligned}
  \|\partial_x^\alpha(uv)\|
  & \lesssim \|u\|_{s}\| v\|_s, \\
  \|[\partial_x^\alpha, u]v\|
  & \lesssim \|\nabla u\|_{s-1}\| v\|_{s-1}.
\end{aligned}
\end{equation*}
\end{lemma}
\begin{lemma}[\cite{H,JJL}]\label{composition-x}
Let $s\in \mathbb{N}$.
Assume that $f(u)$ is a smooth function and $u\in H^s(\mathbb{R}^3)\cap L^{\infty}(\mathbb{R}^3)$.
Then $f(u)\in H^s(\mathbb{R}^3)$, and for any $\alpha$ satisfying $1\leq|\alpha|\leq s$,
\begin{equation*}
  \|\partial_x^\alpha f(u)\|
  \lesssim |\nabla_u f|_{s-1}\|u\|_{L^{\infty}}^{s-1}\|u\|_s.
\end{equation*}
Here $|\cdot|_r$ is the $C^r-$norm.
\end{lemma}

By a straightforward calculation, the above results can be generalized in weighted function spaces as follows.
\begin{corol}\label{Moser-tx}
Let $s>\frac{5}{2}$ be an integer.
Assume that $\|u\|_{s,\e}$ and $\|v\|_{s,\e}$ are bounded.
Then for any $k$ and $\alpha$ satisfying $1\leq k+|\alpha|\leq s$, we have
\begin{equation*}
\begin{aligned}
  \|D^{k,\alpha}(uv)\|& \lesssim \|u\|_{s,\e}\| v\|_{s,\e}, \\
  \|[D^{k,\alpha},u]v\|& \lesssim\|D^1u\|_{s-1,\e}\| v\|_{s-1,\e}.
\end{aligned}
\end{equation*}
\end{corol}

\begin{corol} \label{composition-tx}
Let $s>\frac{3}{2}$ be an integer.
Assume that $f(u)$ is a smooth function  and $\|u\|_{s,\e}$ is bounded.
Then for any $k$ and $\alpha$ satisfying $1\leq k+|\alpha|\leq s$, it holds
\begin{equation*}
\|D^{k,\alpha} f(u)\|
  \lesssim |\nabla_u f|_{s-1}\|u\|_{s,\e}^{s}.
\end{equation*}
\end{corol}

Next, we recall the result on estimating the gradient of vector fields via $\nabla \cdot$ and $\nabla \times$ operators.
\begin{lemma}\label{div-curl}
Let $s \geq 1$ be an integer.
Assume that the smooth vector function $v\in H^s(\mathbb{R}^3)$.
It holds that
\begin{equation*}
  \|\nabla v\|_{s-1}\leq \|\nabla\cdot v\|_{s-1}+\|\nabla \times v\|_{s-1}.
\end{equation*}
\end{lemma}
A simple and direct proof of this result can be found in \cite{LZ}.
When we consider that $v$ lies in a bounded domain with some boundary conditions, the similar results still hold with additional low order and boundary terms. Interested reader can refer to \cite{BB,W} for example.

Now, we show the local existence of the Cauchy problem \eqref{nondim} and \eqref{initial data} for any fixed $\epsilon\in (0,1]$.
We begin with \eqref{guodu}.
Since \eqref{guodu} is a symmetrizable hyperbolic--parabolic system, by selecting the appropriate symmetrizer and following the proof in \cite{K}, we can establish the local existence of the solutions as follows.

\begin{theorem}\label{existence theorem}
Let $s \geq 4$ be an integer.
Assume that the initial data $(\rho_0, u_0,\Theta_0,I_{00},I_{10})$ satisfy
\begin{equation*}
  \|(\rho_0-\underline{\rho}, u_0,\Theta_0-\underline{\Theta},
  I_{00}-\underline{I_0},I_{10}) \|_{s}
  \leq M,
\end{equation*}
for some positive constants $\underline{\rho}$, $\underline{\Theta}$, $\underline{I_0}$ and $M$.
Then there exists a $T^\e>0$ such that the Cauchy problem \eqref{guodu} with the above initial data has a unique classical solution $(\rho,u,\Theta,I_0,I_{1})$ satisfying
$(\rho,u,\Theta,I_0,I_{1})\in C([0,T^\e];H^s(\mathbb{R}^3))$ and
$(\nabla u,\nabla\Theta)\in L^2([0,T^\e];H^{s}(\mathbb{R}^3))$.
Moreover, there exist positive constants $\rho_1$, $\rho_2$, $\Theta_1$, $\Theta_2$, $P_1$ and $P_2$  such that
\begin{equation*}
  \rho_1 \leq \|\rho\|_{L^{\infty}}\leq \rho_2,\;\;
  \Theta_1 \leq \|\Theta\|_{L^{\infty}}\leq \Theta_2,\;\;
  \mathrm{and}\;\;
  P_1 \leq \|P\|_{L^{\infty}}\leq P_2.
\end{equation*}
\end{theorem}
Then, it follows from the transforms \eqref{ptheta} and Theorem \ref{existence theorem} that, for any fixed $\e$ and initial data satisfying \eqref{initial condition}, the Cauchy problem \eqref{nondim} and \eqref{initial data} has a unique classical
solution $(p^\e,u^\e,\theta^{\e},I^\e_0,I^\e_1)$ satisfying
$(p^\e,u^\e,\theta^{\e},I^\e_0,I^\e_1) \in C([0,T^\e];H^s(\mathbb{R}^3))$ and $(\nabla u^\e,\nabla\theta^\e)\in L^2([0,T^\e];H^{s}(\mathbb{R}^3))$.


Finally, we give the local energy decay on the acoustic wave equations obtained by M\'{e}tivier and Schochet \cite{MS} and reformulated by Alazard \cite{A}, which play a important role in acquiring the convergence of solutions.
\begin{lemma}[\cite{A,JJLX}]\label{dispersive}
Let $s>\frac{5}{2}$.
Assume that $v^{\epsilon}\in C([0,T];H^2(\mathbb{R}^3))$  solves the following acoustic wave equation
\begin{equation*}
  \epsilon^2\partial_t(a^{\epsilon}\partial_t v^{\epsilon})-\nabla \cdot (b^{\epsilon}\nabla v^{\epsilon})=c^{\epsilon},
\end{equation*}
where the source term $c^{\epsilon}$ converges to 0 in $L^2([0,T];L^2(\mathbb{R}^3))$ as $\e \rightarrow 0$.
Assume further that the coefficients $(a^{\epsilon},b^{\epsilon})$ are uniformly bounded in $C([0,T];H^{s}(\mathbb{R}^3))$ and converge in $C([0,T];H^s_{\mathrm{loc}}(\mathbb{R}^3))$ to $(a,b)$ as $\e \rightarrow 0$, where $(a,b)$ satisfy the following decay estimates
\begin{align*}
    &|a(x,t)-\underline{a}|\leq N_0|x|^{-1-\sigma},\;\;\;\;|\nabla a(x,t)|\leq N_0|x|^{-2-\sigma},\\
    &|b(x,t)-\underline{b}|\leq N_0|x|^{-1-\sigma},\;\;\;\;|\nabla b(x,t)|\leq N_0|x|^{-2-\sigma},
\end{align*}
for some positive constants $\underline{a},\;\underline{b},\;N_0$ and $\sigma$. Then  $v^{\epsilon}$ converges to 0 in $L^2([0,T];L_{\mathrm{loc}}^2(\mathbb{R}^3))$ as $\e \rightarrow 0$.
\end{lemma}
\section{Uniform estimates}\label{S.3}

According to the local existence of the Cauchy problem \eqref{nondim} and \eqref{initial data} in the previous section, we know that $T^\e$ depends on $\e$, and may tend to zero as $\e \rightarrow 0$.
Using the same argument as in \cite{MS}, to show $T^\e$ has a uniformly positive lower bound, in other words, to establish the uniform-in-time existence of the solutions (Theorem \ref{uniform existence theorem}), it is sufficient to obtain the following a priori estimates.

\begin{prop}\label{priori}
Let $s\geq 4$ be an integer and $(p^\e,u^\e,\theta^\e,I^\e_0,I^\e_1)$ be the classical solution to the Cauchy problem \eqref{nondim} and \eqref{initial data}.
Then there exist constants $\hat{T},\;\e_0\in(0,1]$, and positive increasing polynomial functions $C_0(\cdot)$ and $C(\cdot)$, such that for all $T\in (0,\hat{T}]$ and $\e \in (0,\e_0]$, it holds that
\begin{equation*}
  \mathcal{M}(T)\leq C_0(M_0)+(\sqrt[4]{T}+\e)C(\mathcal{M}(T)),
\end{equation*}
where $\mathcal{M}(T)$ is defined as
\begin{equation*}
  \mathcal{M}(T)=\sup\limits_{t\in[0,T]}\mathcal{M}_1(t)
  +\bigg(\int_0^T \mathcal{M}^2_2(t)\mathrm{d}t \bigg)^{\frac{1}{2}},
\end{equation*}
with
\begin{equation*}
  \mathcal{M}_1(t)=\|(p^\e,u^\e)(t)\|_{s,\e}
     +\|(\e p^\e,\e u^\e, \theta^\e-\theta_c)(t)\|_{s+1,\e}
     +\interleave(I^\e_0-I_c,I^\e_1)(t)\interleave_{s+1,\e},
\end{equation*}
and
\begin{equation*}
  \mathcal{M}_2(t)=\|\nabla(p^\e,u^\e)(t)\|_{s,\e}
      +\|\nabla(\e u^\e,\theta^\e)(t)\|_{s+1,\e}.
\end{equation*}.
\end{prop}

The remainder of this section is devoted to proving Proposition \ref{priori}.
For the sake of notation simplicity, we will drop the superscript ``$\e$" of the variables $(p^\e,u^\e,\theta^{\e},I^\e_0,I^\e_1)$ in the rest of this section.

We first give the following estimates, which can be derived directly from \eqref{nondim}.
\begin{lemma}\label{es-put}
Let $\delta\in(0,2]$ and $s\geq 4$ be an integer. Suppose that $(p,u,\theta,I_0,I_1)$ is a solution to the Cauchy problem \eqref{nondim} and \eqref{initial data} on $[0,T_1]$.
Then, for any $\e \in(0,1]$ and $T\in (0,T_1]$, we have
\begin{equation}\label{deltaI}
  \bigg(\e+\frac{1}{\e^{\delta-1}}\bigg)\|I_1\|_{s,\e}
  \leq \bigg(\e+\frac{1}{\e^{\delta-1}}\bigg)\interleave I_1\interleave_{s,\e}
  \leq C(\mathcal{M}_1),
\end{equation}
and
\begin{equation}\label{es-p u theta t}
  \|(\e \partial_t p,\e \partial_t u,\partial_t \theta )\|_{s,\e}\lesssim C(\mathcal{M}_1)(1+\mathcal{M}_2)+\mathcal{M}_2.
\end{equation}
\end{lemma}

\begin{proof}

We note that if $k=0$, then $D^{k,\alpha}=\dot{D}^{k,\alpha}$;
and if $k>0$, then $D^{k,\alpha}=\e \dot{D}^{k,\alpha}$.
Thus, for any $\e \in(0,1]$, we have
\begin{equation}\label{bijiao}
  \|\cdot\|_{i,\e}\leq \interleave \cdot \interleave_{i,\e},\;\;\;\; i=1,\,2,\, \dots ,\,s+1.
\end{equation}
We rewrite \eqref{nondim}$_5$ as
\begin{equation*}
\bigg(\e+\frac{1}{\e^{\delta-1}}\bigg)I^\epsilon_1=
  -\e \partial_t I^\epsilon_1-\frac{\nabla I^\epsilon_0}{3},
\end{equation*}
and then obtain
\begin{equation}\label{11}
  \bigg(\e+\frac{1}{\e^{\delta-1}}\bigg)\interleave I_1\interleave_{s,\e}
  \leq C(\mathcal{M}_1).
\end{equation}
Combining \eqref{bijiao} and \eqref{11}, we get \eqref{deltaI}.

According to the equation of $p$ in \eqref{nondim}, it follows that
\begin{align*}
   \|\e \partial_t p\|_{s,\e}
    \lesssim{}& \|\e u\|_{s,\e}\| \nabla p\|_{s,\e}+\|\nabla u\|_{s,\e}
    +\|(\e p,\theta-\theta_c)\|^{s}_{s,\e}
    \big(\|\nabla \theta\|^{2}_{s,\e}+\|\nabla \theta\|_{s+1,\e}\big)\\
    &  +\|\e p\|^{s}_{s,\e}\big[\|\e u\|^2_{s+1,\e}
        + \|I_0-I_c\|_{s,\e}+ \|\theta-\theta_c\|^s_{s,\e}
       +(\e^2+\e^{2-\delta})\|I_1\|_{s,\e}\| u\|_{s,\e}\big]\\
    \lesssim{}& C(\mathcal{M}_1)\big(1+\| \nabla p\|_{s,\e}
            +\|\nabla \theta\|_{s+1,\e}\big)+\|\nabla u\|_{s,\e}\\
    \leq {}&C(\mathcal{M}_1)(1+\mathcal{M}_2)+\mathcal{M}_2.
\end{align*}
The estimates of $ \e\partial_t u $ and $\partial_t \theta$ can be handled in the same way and we omit the details for simplicity.
\end{proof}

Now, we give the $s+1$ order derivative estimates of $(\e p,\e u, \theta-\theta_c,I_0-I_c,I_1)$.

\begin{lemma}\label{ep-eu-I}
Let  $\delta\in (0,2]$ and $s\geq 4$ be an integer.
Suppose that $(p,u,\theta,I_0,I_1)$ is a solution to the Cauchy problem \eqref{nondim} and \eqref{initial data} on $[0,T_1]$.
Then for any $\e \in(0,1]$ and $T\in (0,\hat{T}]$, $\hat{T}=\mathrm{min}\{T_1,1\}$, the following estimates hold.
\begin{equation}\label{I0I1}
\interleave(I_0-I_c,I_1)\interleave_{s+1,\epsilon,T}
   +\bigg(\int_0^T\interleave(I_0-I_c)(t)\interleave^{2}_{s+1,\epsilon}
   + \bigg(1+\frac{1}{\e^\delta}\bigg)\interleave I_1(t)\interleave^{2}_{s+1,\epsilon}
   \mathrm{d}t \bigg)^{\frac{1}{2}}
   \leq C_0(M_0)+\sqrt[4]{T}C(\mathcal{M}),
\end{equation}
\begin{equation}\label{hat-p-u-theta}
   \|(\e p,\e u,\theta-\theta_c)\|_{s+1,\e,T}
   +\bigg(\int_0^T\|\nabla(\e u, \theta)(t)\|_{s+1,\e}^2 \mathrm{d}t\bigg)^{\frac{1}{2}}
   \leq C_0(M_0)+\sqrt[4]{T}C(\mathcal{M}).
\end{equation}
\end{lemma}

\begin{proof}
We first construct auxiliary equations to obtain the anti-symmetric structure, which help us to cancel the singular terms, and to reduce the order of the singularity of $\frac{1}{\e^{\delta}}I_1$ in \eqref{nondim}$_2$.
Setting
\begin{equation*}
  (\hat{p},\hat{u},\hat{\theta},\hat{I}_0,\hat{I}_1)=(\e p-(\theta-\theta_c),\e u,\theta-\theta_c,I_0-I_c,I_1),
\end{equation*}
which means that $(p, u,\theta,I_0,I_1)=(\frac{\hat{p}+\hat{\theta}}{\e}, \frac{\hat{u}}{\e},\hat{\theta}+\theta_c,\hat{I}_0+I_c,\hat{I}_1)$, and putting them into \eqref{nondim}, we get from a straightforward calculation that $(\hat{p},\hat{u},\hat{\theta},\hat{I}_0,\hat{I}_1)$ solves the following auxiliary equations
 \begin{equation}\label{nondim-hat}
\left\{
\begin{aligned}
  &  \partial_t \hat{p}+u \cdot \nabla \hat{p}+\frac{1}{\epsilon}\nabla \cdot \hat{u}= 0, \\
  &  e^{-\theta} (\partial_t \hat{u} + u\cdot \nabla \hat{u} )
       + \frac{1}{\epsilon}(\nabla \hat{p}+\nabla \hat{\theta})
     =e^{-\e p}\bigg[\nabla\cdot \Psi(\hat{u})+\bigg(\e+\frac{1}{\e^{\delta-1}}\bigg) \hat{I}_1\bigg],\\
  & \partial_t \hat{\theta} + u\cdot \nabla \hat{\theta}+ \frac{1}{\e} \nabla \cdot \hat{u}
    = e^{- \e p}\big[\kappa e^{\theta_c}\Delta e^{\hat{\theta}}
       +\Psi(\hat{u}):\nabla \hat{u}
       +\epsilon\hat{I}_0-\e e^{4\theta_c}(e^{4\hat{\theta}}-1)
       -(\e+\e^{1-\delta}) \hat{I}_1\cdot \hat{u}\big],  \\
  &  \partial_t \hat{I}_0+\frac{1}{\epsilon}\nabla \cdot \hat{I}_1
      =e^{4\theta_c}(e^{4\hat{\theta}}-1)-\hat{I}_0,\\
  & \partial_t \hat{I}_1+\frac{1}{3\epsilon}\nabla \hat{I}_0=-\bigg(1+\frac{1}{\e^\delta}\bigg)\hat{I}_1,
\end{aligned}
\right.
\end{equation}
with the initial data
\begin{equation*}
(\hat{p},\hat{u},\hat{\theta},\hat{I}_{0},\hat{I}_{1})|_{t=0}
  =(\e p_0-(\theta_0-\theta_c),\e u_0,\theta_0-\theta_c,I_{00}-I_c,I_{10}).
\end{equation*}

For the above system, applying $\dot{D}^{k,\alpha}$, $0\leq k+|\alpha|\leq s+1$, to \eqref{nondim-hat}$_4$-\eqref{nondim-hat}$_5$, we get
\begin{equation}\label{DI}
\left\{
\begin{aligned}
  &  \partial_t \dot{D}^{k,\alpha}\hat{I}_0
       +\frac{1}{\epsilon}\nabla \cdot \dot{D}^{k,\alpha} \hat{I}_1
    =e^{4\theta_c}\dot{D}^{k,\alpha}(e^{4\hat{\theta}}-1)-\dot{D}^{k,\alpha}\hat{I}_0,\\
  & \partial_t \dot{D}^{k,\alpha}\hat{I}_1+\frac{1}{3\e}\nabla \dot{D}^{k,\alpha}\hat{I}_0
    =-\bigg(1+\frac{1}{\e^{\delta}}\bigg)\dot{D}^{k,\alpha}\hat{I}_1.
\end{aligned}
\right.
\end{equation}
Multiplying \eqref{DI} with $(2\dot{D}^{k,\alpha} \hat{I}_0, 6\dot{D}^{k,\alpha} \hat{I}_1)$, integrating the result over $\mathbb{R}^3$ and using integration by parts, we obtain that
\begin{equation}\label{dt-I0I1}
   \frac{\mathrm{d}}{\mathrm{d}t}
     \big(\|\dot{D}^{k,\alpha}\hat{I}_0\|^2
         +3\|\dot{D}^{k,\alpha}\hat{I}_1\|^2\big)
         +2\|\dot{D}^{k,\alpha}\hat{I}_0\|^2
         +6\bigg(1+\frac{1}{\e^\delta}\bigg)\|\dot{D}^{k,\alpha}\hat{I}_1\|^2
  \lesssim \|\dot{D}^{k,\alpha}\hat{I}_0\|\| \dot{D}^{k,\alpha}(e^{4\hat{\theta}}-1)\|.
\end{equation}
Thanks to the estimate of $\partial_t \theta$ in \eqref{es-p u theta t}, it follows that
\begin{align}
    \|\dot{D}^{k,\alpha}\hat{I}_0\|\|\dot{D}^{k,\alpha}(e^{4\hat{\theta}}-1)\|
  \leq {}& \|\dot{D}^{k,\alpha}\hat{I}_0\|
      \big(\|(e^{4\hat{\theta}}-1)\|_{s+1}
      +\|e^{4\hat{\theta}}\partial_t \hat{\theta}\|_{s,\e}\big)\notag\\
  \lesssim {}& \interleave I_0-I_c\interleave_{s+1,\e}\big(\|\theta-\theta_c\|^{s+1}_{s+1,\e}
      +\|\theta-\theta_c\|^{s}_{s,\e}\|\partial_t \theta\|_{s,\e}\big)\notag\\
  \leq {}& C(\mathcal{M}_1)(1+\mathcal{M}_2).             \label{1}
\end{align}
Substituting \eqref{1} into \eqref{dt-I0I1}, we get
\begin{equation}\label{eI1}
   \frac{\mathrm{d}}{\mathrm{d}t}
     \big(\|\dot{D}^{k,\alpha}\hat{I}_0\|^2
         +3\|\dot{D}^{k,\alpha}\hat{I}_1\|^2\big)
         +2\|\dot{D}^{k,\alpha}\hat{I}_0\|^2
         +6\bigg(1+\frac{1}{\e^\delta}\bigg)\|\dot{D}^{k,\alpha}\hat{I}_1\|^2
  \leq C(\mathcal{M}_1)(1+\mathcal{M}_2).
\end{equation}
The integration of \eqref{eI1} on $[0,T]$ yields
\begin{align*}
   & \|\dot{D}^{k,\alpha}\hat{I}_0\|^2+3\|\dot{D}^{k,\alpha}\hat{I}_1\|^2
         +\int_0^T 2\|\dot{D}^{k,\alpha}\hat{I}_0\|^2
         +6\bigg(1+\frac{1}{\e^\delta}\bigg)\|\dot{D}^{k,\alpha}\hat{I}_1\|^2
         \mathrm{d}t \notag\\
  \leq {}& \|\dot{D}^{k,\alpha}I_{00}-I_c\|^2+3\|\dot{D}^{k,\alpha}I_{10}\|^2
                +TC(\mathcal{M}_1)+\sqrt{T}C(\mathcal{M}_1)
         \bigg(\int_0^T\mathcal{M}_2^2\mathrm{d}t\bigg)^{\frac{1}{2}} \\
   \leq {}& C_0(M_0)+\sqrt{T}C(\mathcal{M}).
\end{align*}
Summing up the above result for all $k$ and $\alpha$ satisfying $0\leq k+|\alpha|\leq s+1$, we obtain the desired estimate \eqref{I0I1}.

Similarly, applying $D^{k,\alpha}$, $0\leq k+|\alpha|\leq s+1$, to \eqref{nondim-hat}$_1$-\eqref{nondim-hat}$_3$, we have
\begin{equation}\label{Dka}
\left\{
\begin{aligned}
  &  \partial_t D^{k,\alpha}\hat{p}+u \cdot \nabla D^{k,\alpha}\hat{p}
         +\frac{1}{\epsilon}\nabla \cdot D^{k,\alpha}\hat{u}= g_1,\\
  &  e^{-\theta} (\partial_t D^{k,\alpha}\hat{u} + u\cdot \nabla D^{k,\alpha}\hat{u} )
       + \frac{1}{\epsilon}(\nabla D^{k,\alpha}\hat{p}+\nabla D^{k,\alpha}\hat{\theta})
     =e^{-\e p}\nabla\cdot \Psi(D^{k,\alpha}\hat{u})
        +g_2,\\
  & \partial_t D^{k,\alpha}\hat{\theta} + u\cdot \nabla D^{k,\alpha}\hat{\theta}
     + \frac{1}{\e}\nabla \cdot D^{k,\alpha}\hat{u}
    = \kappa e^{-\e p+\theta}\Delta D^{k,\alpha}\hat{\theta}  +g_3,
\end{aligned}
\right.
\end{equation}
where
\begin{align*}
g_1= & -\big[D^{k,\alpha},u\cdot \nabla\big] \hat{p},\\
g_2= & -\big [D^{k,\alpha}, e^{-\theta}(\partial_t+ u\cdot \nabla)\big]\hat{u}
         +\big [D^{k,\alpha},e^{-\e p} \big]\nabla\cdot\Psi(\hat{u})
         +\bigg(\e+\frac{1}{\e^{\delta-1}}\bigg) D^{k,\alpha}(e^{-\e p}\hat{I}_1),\\
g_3= & -\big [D^{k,\alpha}, u\cdot \nabla \big]\hat{\theta}
         +\big [D^{k,\alpha},\kappa e^{-\e p+\theta}\big] \Delta \hat{\theta}
         +D^{k,\alpha}(\kappa e^{-\e p+\theta} \nabla\hat{\theta}\cdot\nabla\hat{\theta})\\
     & + D^{k,\alpha}\big\{ e^{-\e p}\big[\Psi(\hat{u}):\nabla \hat{u}+\epsilon\hat{I}_0
                 -\e e^{4\theta_c}(e^{4\hat{\theta}}-1)
                 - (\e+\e^{1-\delta}) \hat{I}_1\cdot \hat{u} \big] \big\}.
\end{align*}
Multiplying \eqref{Dka} with
$(2D^{k,\alpha} \hat{p}, 2D^{k,\alpha} \hat{u}, 2D^{k,\alpha} \hat{\theta})$,
integrating over $\mathbb{R}^3$ and using integration by parts, we arrive at

\begin{equation}\label{dt-p-u-theta}
   \frac{\mathrm{d}}{\mathrm{d}t}\int_{\mathbb{R}^3}
       |D^{k,\alpha}\hat{p}|^2
       +e^{-\theta}|D^{k,\alpha}\hat{u}|^2
       +|D^{k,\alpha}\hat{\theta}|^2\mathrm{d}x
     \leq \sum\limits_{i=1}^5 J_i,
\end{equation}
where
\begin{align*}
   J_1={}& \int_{\mathbb{R}^3}(\partial_t e^{-\theta})|D^{k,\alpha}\hat{u}|^2 \mathrm{d}x,\\
   J_2={}& \int_{\mathbb{R}^3}(\nabla \cdot u)|D^{k,\alpha}\hat{p}|^2
           +\nabla \cdot (e^{-\theta}u)|D^{k,\alpha}\hat{u}|^2
           +(\nabla \cdot u)|D^{k,\alpha}\hat{\theta}|^2\mathrm{d}x,\\
   J_3={}& 2\big\langle D^{k,\alpha}\hat{u},
           e^{-\e p}\nabla\cdot \Psi(D^{k,\alpha}\hat{u})\big\rangle,\\
   J_4={}& 2\big\langle D^{k,\alpha}\hat{\theta},
          \kappa e^{-\e p+\theta}\Delta D^{k,\alpha}\hat{\theta}
          \big\rangle, \\
   J_5={}& 2\big\langle D^{k,\alpha}\hat{p},g_1\big\rangle
           +\big\langle D^{k,\alpha}\hat{u},g_2\big\rangle
           +\big\langle D^{k,\alpha}\hat{\theta},g_3\big\rangle.
\end{align*}

Now, we estimate these terms in turn.
Using the equation of $\theta$ in \eqref{nondim}, it is easy to see that
$$\|\partial_t \theta\|_2 \leq C(\mathcal{M}_1).$$
Then, we get from Sobolev's inequality that
\begin{equation}\label{j1}
  |J_1|= \left|\int_{\mathbb{R}^3}(\partial_t e^{-\theta})|D^{k,\alpha}\hat{u}|^2 \mathrm{d}x\right|
  \leq  \|e^{-\theta}\|_{L^{\infty}}
  \|\partial_t\theta\|_{L^{\infty}}
  \|D^{k,\alpha}\hat{u}\|^2
  \leq C(\mathcal{M}_1).
\end{equation}
For the term $J_2$,  we have
\begin{align}
  |J_2| \leq {} & \|\nabla \cdot u\|_{L^{\infty}}\|D^{k,\alpha}\hat{p}\|^2
            +\|\nabla \cdot (e^{-\theta}u)\|_{L^{\infty}}\|D^{k,\alpha}\hat{u}\|^2
            +\|\nabla \cdot u\|_{L^{\infty}}\|D^{k,\alpha}\hat{\theta}\|^2\notag\\
  \lesssim {} & \|u\|_3\|D^{k,\alpha}\hat{p}\|^2
            +(\|\nabla\theta\|_2\|u\|_2+\|u\|_3)\|e^{-\theta}\|_{L^{\infty}}\|D^{k,\alpha}\hat{u}\|^2
            +\|u\|_3\|D^{k,\alpha}\hat{\theta}\|^2\notag\\
  \leq {} & C(\mathcal{M}_1). \label{j2}
\end{align}
Using  integration by parts, we obtain that
\begin{align}
  J_3={}& 2\big\langle D^{k,\alpha}\hat{u},
           e^{-\e p}[\mu \Delta D^{k,\alpha}\hat{u}
           +(\mu+\lambda)\nabla \nabla \cdot D^{k,\alpha}\hat{u}]\big\rangle \notag\\
     ={}& -2\big\langle \nabla D^{k,\alpha}\hat{u},
             \mu e^{-\e p}\nabla D^{k,\alpha}\hat{u}\big\rangle
          -2\big\langle \nabla\cdot D^{k,\alpha}\hat{u},
           (\mu+\lambda)e^{-\e p} \nabla \cdot D^{k,\alpha}\hat{u}\big\rangle \notag\\
        & -2\big\langle D^{k,\alpha}\hat{u},
           \mu\nabla e^{-\e p}\cdot \nabla D^{k,\alpha}\hat{u}
           +(\mu+\lambda)\nabla e^{-\e p}\cdot
            \nabla  D^{k,\alpha}\hat{u}\big\rangle \notag\\
  \leq {}& -l_1\| D^{k,\alpha}\nabla\hat{u}\|^2
           +(2\mu+\lambda)\|\nabla e^{-\e p}\|_{L^{\infty}} \| D^{k,\alpha}\hat{u}\|
              \| D^{k,\alpha}\nabla\hat{u}\| \notag\\
  \leq {}& -l_1\| D^{k,\alpha}\nabla\hat{u}\|^2
              +C(\mathcal{M}_1)\mathcal{M}_2,    \label{j3}
\end{align}
for some constant $l_1>0$.
Similarly, for the term $J_4$, we have
\begin{align}
J_4={}& 2\big\langle D^{k,\alpha}\hat{\theta},
          \kappa e^{-\e p+\theta}\Delta D^{k,\alpha}\hat{\theta}\big\rangle \notag\\
   ={}& -2\big\langle \nabla D^{k,\alpha}\hat{\theta},
            \kappa e^{-\e p+\theta} \nabla D^{k,\alpha}\hat{\theta}\big\rangle
        -2\big\langle D^{k,\alpha}\hat{\theta},
          \kappa \nabla e^{-\e p+\theta}\cdot
           D^{k,\alpha} \nabla\hat{\theta}\big\rangle  \notag\\
   \leq{}&  -l_2\| D^{k,\alpha}\nabla\hat{\theta}\|^2
            +C(\mathcal{M}_1)\mathcal{M}_2, \label{j4}
\end{align}
for some  positive constant $l_2$.

It remains to estimate the term $ J_5$.
According to  Corollary \ref{Moser-tx} and Lemma \ref{es-put}, it follows that
\begin{align*}
 \|D^{k,\alpha}\hat{p}\|\|g_1\|
  = {}&\|D^{k,\alpha}\hat{p}\|\big\|\big[D^{k,\alpha},u\cdot \nabla\big] \hat{p}\big\| \\
  \lesssim{}& \|\hat{p}\|_{s+1,\e}\|D^1 u\|_{s,\e}\|\hat{p}\|_{s+1,\e}\\
  \lesssim{}& \|\hat{p}\|^2_{s+1,\e}(\|\e\partial_t u\|_{s,\e}+\|\nabla u\|_{s,\e})\\
  \leq{}& C(\mathcal{M}_1)(1+\mathcal{M}_2),
\end{align*}
and
\begin{align*}
   \|D^{k,\alpha}\hat{u}\|\|g_2\|
   \leq{}& \|D^{k,\alpha}\hat{u}\|
             \big\|\big[D^{k,\alpha}, e^{-\theta}(\partial_t+ u\cdot \nabla)\big]\hat{u}\big\|
         + \|D^{k,\alpha}\hat{u}\|
             \big\|\big [D^{k,\alpha},e^{-\e p}\big]\nabla\cdot\Psi(\hat{u})\big\|\\
        {}& + \bigg(\e+\frac{1}{\e^{\delta-1}}\bigg)\|D^{k,\alpha}\hat{u}\|
              \big\| D^{k,\alpha}\big(e^{-\e p}\hat{I}_1\big)\big\|\\
  \lesssim {}& \|\hat{u}\|_{s+1,\e}\|D^1e^{-\theta}\|_{s,\e}\|\partial_t \hat{ u}\|_{s,\e}
              +\|\hat{u}\|_{s+1,\e}\|D^1(e^{-\theta}u)\|_{s,\e}\|\nabla\hat{u}\|_{s,\e}\\
           {}&+\|\hat{u}\|_{s+1,\e}\|D^1e^{-\e p}\|_{s,\e}\|\nabla \hat{u}\|_{s+1,\e}
              +\bigg(\e+\frac{1}{\e^{\delta-1}}\bigg)\|\hat{u}\|_{s+1,\e}
                \|\e p\|^{s+1}_{s+1,\e}\|\hat{I}_1\|_{s+1,\e}\\
  \lesssim {}& \|\hat{u}\|_{s+1,\e}\|\theta-\theta_c\|^{s+1}_{s+1,\e}
              \big[\|\e \partial_t u\|_{s,\e}+(\|\e \partial_t u\|_{s,\e}
                    +\|\nabla u\|_{s,\e})\|\nabla\hat{u}\|_{s,\e}\big]\\
           {}&+ \|\hat{u}\|_{s+1,\e}\|\e p\|^{s+1}_{s+1,\e}\|\nabla \hat{u}\|_{s+1,\e}
              + \|\hat{u}\|^2_{s+1,\e}\|\e p\|^{2s+2}_{s+1,\e}
              +\bigg(\e+\frac{1}{\e^{\delta-1}}\bigg)^2\|\hat{I}_1\|^2_{s+1,\e}\\
   \leq {}& C(\mathcal{M}_1)(1+\mathcal{M}_2)
               +\frac{1}{\e^{2\delta-2}}\interleave\hat{I}_1\interleave^{2}_{s+1,\e}.
\end{align*}
In the same way, we have
\begin{equation*}
\| D^{k,\alpha}\hat{\theta}\|\|g_3\|
\leq  C(\mathcal{M}_1)(1+\mathcal{M}_2)
+\frac{1}{\e^{2\delta-2}}\interleave\hat{I}_1\interleave^{2}_{s+1,\e}.
\end{equation*}
Then, we obtain that
\begin{align}
   |J_5|\leq{}& 2\| D^{k,\alpha}\hat{p}\|\|g_1\|
           +2\| D^{k,\alpha}\hat{u}\|\|g_2\|
           +2\| D^{k,\alpha}\hat{\theta}\|\|g_3\|\notag\\
   \leq {}& C(\mathcal{M}_1)(1+\mathcal{M}_2)
           +\frac{4}{\e^{2\delta-2}}\interleave \hat{I}_1\interleave^{2}_{s+1,\e}. \label{j5}
\end{align}
Putting \eqref{j1}-\eqref{j4} and \eqref{j5} into \eqref{dt-p-u-theta}, and integrating the resulting inequality on $[0,T]$ gives
\begin{align}
&\|D^{k,\alpha} \hat{p} \|+\|e^{-\frac{\theta}{2}}D^{k,\alpha} \hat{u} \|+\|D^{k,\alpha} \hat{\theta} \|
+\int_0^T l_1\|D^{k,\alpha} \nabla \hat{u} \| +l_2\|D^{k,\alpha} \nabla \hat{\theta} \| \mathrm{d}t \notag\\
\leq {}&C_0(\mathcal{M}_0)+\sqrt{T}C(\mathcal{M})
+\int_0^T\frac{4}{\e^{2\delta-2}}\interleave\hat{I}_1\interleave^{2}_{s+1,\e} \mathrm{d}t. \label{C}
\end{align}
Since $\delta\in(0,2]$, we have for any $\e\in(0,1]$ that
\begin{equation*}
  \frac{1}{\e^{2\delta-2}}\leq \frac{1}{\e^{\delta}}.
\end{equation*}
Recalling the boundedness of $\int_0^T \frac{1}{\e^\delta}\interleave \hat{I}_1\interleave^{2}_{s+1,\e}\mathrm{d}t$ in \eqref{I0I1}
and then summing up \eqref{C} for all $k$ and $\alpha$ satisfying $0\leq k+|\alpha|\leq s+1$, we obtain \eqref{hat-p-u-theta}.
\end{proof}

Next, we give the estimate on $\|(\e\partial_t)^k(p,u)\|$, which is a part of the weighted norm, and it also helps us get the estimates on $(\nabla p, \nabla \cdot u)$.

\begin{lemma}\label{pt-ut-lemma}
Let $\delta\in (0,2]$ and $s\geq 4$ be an integer.
Assume that $(p,u,\theta,I_0,I_1)$ is a solution to the Cauchy problem \eqref{nondim} and \eqref{initial data} on $[0,T_1]$.
Then for any $\e \in(0,1]$ and $T\in (0,\hat{T}]$, $\hat{T}=\mathrm{min}\{T_1,1\}$, it holds that
\begin{equation}\label{pk-uk}
   \sup \limits_{t\in[0,T]}\|(\e\partial_t)^k( p, u)\|
  +\bigg(\int_0^T\| \nabla (\e\partial_t)^k u(t)\|^2\mathrm{d}t\bigg)^{\frac{1}{2}}
    \leq C_0(M_0)+\sqrt[4]{T}C(\mathcal{M}),\;\; k=0,1, \dots , s.
\end{equation}
\end{lemma}

\begin{proof}

According to the singular structure of \eqref{nondim}$_1$-\eqref{nondim}$_2$, we introduce the equivalent pressure and velocity:
\begin{equation}\label{bianhuan}
  \tilde{p}=p+\frac{e^{-\e p}(I_0-I_c)}{3}\;\;\mathrm{and}\;\;\tilde{u}=2u-\kappa e^{-\epsilon p+\theta}\nabla\theta.
\end{equation}
With the help of \eqref{nondim}$_3$ and \eqref{nondim}$_5$, we construct the equations of
$(\tilde{p},\tilde{u})$ in the following form, which own an anti-symmetric structure,
\begin{equation}\label{pt-ut}
\left\{
\begin{aligned}
  &  \partial_t \tilde{p}+u\cdot \nabla \tilde{p}
     +\frac{1}{\epsilon}\nabla \cdot\tilde{u}
   = g_4, \\
  &  \frac{e^{-\theta}}{2} (\partial_t \tilde{u} + u\cdot \nabla \tilde{u} )
       + \frac{\nabla \tilde{p}}{\epsilon}
     =\frac{e^{-\e p}}{2}\Big(\nabla\cdot \Psi(\tilde{u})
        +\frac{\kappa}{2}\nabla \nabla \cdot\tilde{u}\Big) +g_5,
\end{aligned}
\right.
\end{equation}
where
\begin{align*}
  g_4={} & \frac{e^{-\e p}}{3}\big[-\e (I_0-I_c)(\partial_t p+ u\cdot \nabla p)
     +(\partial_t I_0+ u\cdot \nabla I_0)\big]\\
  {} & +e^{-\e p}\bigg[\e \Psi(u):\nabla u + (I_0-e^{4\theta})
     -\bigg(\e+\frac{1}{\e^{\delta-1}}\bigg) I_1\cdot u\bigg]
     +\kappa e^{-\e p+\theta}\nabla p\cdot \nabla \theta,
\end{align*}
and
\begin{align*}
  g_5={} &\frac{\kappa}{2}e^{-\e p}\big[\e(\partial_t p+u\cdot \nabla p)\nabla \theta
              -(\partial_t \theta+u\cdot \nabla \theta)\nabla \theta\big]\\
       {}&+\frac{\kappa }{2}e^{-\e p}\bigg[\mu\Delta(e^{-\e p+\theta}\nabla \theta)
         +\bigg(\mu+\lambda+\frac{\kappa}{2}\bigg)\nabla\nabla\cdot( e^{-\e p+\theta}\nabla \theta)
         -\kappa \nabla (e^{-\e p}\Delta e^{\theta})\bigg]\\
       {} &+e^{-\e p}\bigg[\frac{\kappa}{2}\nabla u \nabla \theta
                   -\partial_t I_1-\frac{(I_0-I_c) \nabla p}{3}\bigg]\\
       {} &-\frac{\kappa}{2}e^{-\e p}\nabla \big \{e^{-\e p}\big[\epsilon^2\Psi(u):\nabla u
              +\epsilon(I_0-e^{4\theta})
              -(\epsilon^2+\e^{2-\delta}) I_1\cdot u\big]\big\}.
\end{align*}
Let $(\tilde{p}^k,\tilde{u}^k)=(\e\partial_t)^k(\tilde{p},\tilde{u})$.
Applying the operator $(\e\partial_t)^k$, $k=0,1, \dots ,s$, to \eqref{pt-ut} and taking the $L^2$-inner product of the resulting equations with $2(\tilde{p}^k,\tilde{u}^k)$ yields
\begin{equation}\label{dt-p-u}
   \frac{\mathrm{d}}{\mathrm{d}t}\int_{\mathbb{R}^3}
       |\tilde{p}^k|^2+\frac{e^{-\theta}}{2}|\tilde{u}^k|^2\mathrm{d}x
       \leq \sum\limits_{i=1}^5 K_i,
\end{equation}
where
\begin{align*}
   K_1={}& \frac{1}{2}\int_{\mathbb{R}^3}(\partial_t e^{-\theta})|\tilde{u}^k|^2 \mathrm{d}x,\\
   K_2={}& \int_{\mathbb{R}^3}(\nabla \cdot u)|\tilde{p}^k|^2
           +\frac{1}{2}\nabla \cdot (e^{-\theta}u)|\tilde{u}^k|^2\mathrm{d}x,\\
   K_3={}& \frac{1}{2}\big\langle \tilde{u}^k,e^{-\e p}[2\nabla\cdot \Psi(\tilde{u}^k)
                           +\kappa\nabla \nabla \cdot\tilde{u}]\big\rangle,\\
   K_4={}& \big\langle 2\tilde{p}^k, (\e \partial_t)^k g_4\big\rangle
            +\big\langle 2\tilde{p}^k,-\big[(\e \partial_t)^k,u\cdot \nabla\big]
            \tilde{p}\big\rangle,\\
   K_5= {}& \big\langle 2\tilde{u}^k,(\e \partial_t)^k g_5\big\rangle
              +\big\langle \tilde{u}^k,-\big [(\e \partial_t)^k,
                e^{-\theta}(\partial_t+ u\cdot \nabla)\big]\tilde{u}\big\rangle\\
        {}&    +\big\langle \tilde{u}^k, \big [(\e \partial_t)^k,
               e^{-\e p} \big]\nabla\cdot\Psi(\tilde{u})\big\rangle
               +\frac{\kappa}{2}\big\langle \tilde{u}^k, \big [(\e \partial_t)^k,
               e^{-\e p} \big]\nabla\nabla \cdot\tilde{u}\big\rangle.
\end{align*}
Utilizing an argument similar to the one used to bound the terms $J_i$ $(i=1, \dots ,5)$ in Lemma \ref{ep-eu-I}, we can get
\begin{align}
  |K_1| \leq {}&\|\tilde{u}^k\|^2 C(\mathcal{M}_1),\label{k1}\\
  |K_2| \leq{}&\big(\|p^k\|^2+\|\tilde{u}^k\|^2\big) C(\mathcal{M}_1),\label{k2}\\
  K_3 \leq {}& -l_3\| \nabla\tilde{u}^k\|^2+C(\mathcal{M}_1)(1+\mathcal{M}_2),\label{k3}\\
  |K_4| \leq  {}& C(\mathcal{M}_1)(1+\mathcal{M}_2),\label{k4}
\end{align}
where $l_3>0$ is a constant.

Now, let us focus on $K_5$.
For the first term in $K_5$, we rewrite it as
\begin{equation}\label{g5}
  \big\langle 2\tilde{u}^k,(\e \partial_t)^k g_5\big\rangle
    =\sum\limits_{i=1}^7 Q_i,
\end{equation}
where
\begin{align*}
  Q_1=  {} &\kappa \big\langle \tilde{u}^k, (\e \partial_t)^k
  \big[ \e e^{-\e p}(\partial_t p+u\cdot \nabla p)\nabla \theta\big]\big\rangle, \\
  Q_2=  {} &-\kappa\big\langle \tilde{u}^k, (\e \partial_t)^k
  \big[ e^{-\e p}(\partial_t \theta+u\cdot \nabla \theta)\nabla \theta\big]\big\rangle,\\
  Q_3=  {} &\kappa\mu \big\langle \tilde{u}^k,(\e \partial_t)^k
  \big[ e^{-\e p}\Delta(e^{-\e p+\theta}\nabla \theta)\big]\big\rangle,\\
  Q_4= {} &\frac{1}{2}\kappa(2\mu+2\lambda+\kappa)\big\langle \tilde{u}^k,(\e \partial_t)^k
  \big[e^{-\e p}\nabla\nabla\cdot(e^{-\e p+\theta}
             \nabla \theta)\big]\big\rangle,\\
  Q_5={} &-\kappa^2\big\langle \tilde{u}^k,(\e \partial_t)^k
  \big[e^{-\e p}\nabla(e^{-\e p}\Delta e^{\theta})\big]\big\rangle,\\
  Q_6={} &\frac{1}{3}\big\langle \tilde{u}^k,(\e \partial_t)^k
  [e^{-\e p}(3\kappa\nabla u \nabla \theta -6\partial_t I_1-2(I_0-I_c) \nabla p)]\big\rangle,
\end{align*}
and
\begin{equation*}
  Q_7= -\kappa\big\langle \tilde{u}^k,(\e \partial_t)^k
  \big\{e^{-\e p}\nabla \big[e^{-\e p}\big(\epsilon^2\Psi(u):\nabla u
              +\epsilon(I_0-e^{4\theta})
              -(\epsilon^2+\e^{2-\delta}) I_1\cdot u\big)\big]\big\}\big\rangle.
\end{equation*}
First, according to Lemma \ref{es-put} and the fact that
\begin{equation*}
  \|\tilde{u}^k\|=\|2(\e \partial_t)^ku\|+\|(\e \partial_t)^k(\kappa e^{-\epsilon p+\theta}\nabla\theta)\|\leq C(\mathcal{M}_1),
\end{equation*}
the term $Q_1$ can be bounded as follows
\begin{align*}
   |Q_1|\lesssim {}& \|\tilde{u}^k\|\|e^{-\e p}\nabla \theta\|_{s,\e}
              \|\e \partial_t p+\e u \cdot \nabla p \|_{s,\e}\\
     \lesssim {}& \|\tilde{u}^k\|\|\e p\|^s_{s,\e}\|\theta-\theta_c\|_{s+1,\e}
              (\|\e \partial_t  p\|_{s,\e}+\|u\|_{s,\e}\|\e p\|_{s+1,\e})\\
    \leq {}&C(\mathcal{M}_1)(1+\mathcal{M}_2).
\end{align*}
In the same way, we see that
\begin{align*}
   |Q_2|\lesssim {}& \|\tilde{u}^k\|\|e^{-\e p}\nabla \theta\|_{s,\e}
              \| \partial_t \theta+ u \cdot \nabla \theta \|_{s,\e}\\
    \leq {}&C(\mathcal{M}_1)(1+\mathcal{M}_2).
\end{align*}
Next, since
\begin{equation*}
    \|\nabla\tilde{u}^k\|=\|2(\e \partial_t)^k\nabla u\|+\|(\e \partial_t)^k\nabla
    (\kappa e^{-\epsilon p+\theta}\nabla\theta)\|\leq C(\mathcal{M}_1)(1+\mathcal{M}_2)+\mathcal{M}_2,
\end{equation*}
using integrating by parts and Cauchy's inequality, we get the estimate of $Q_3$ as follows,
\begin{align*}
  Q_3
   ={} &\kappa \mu\big\langle \tilde{u}^k,(\e \partial_t)^k
  \big[ e^{-\e p}\Delta(e^{-\e p+\theta}\nabla \theta)\big]\big\rangle\\
   ={} &-\kappa \mu \big\langle \nabla \tilde{u}^k,(\e \partial_t)^k
  \big[ e^{-\e p}\nabla(e^{-\e p+\theta}\nabla \theta)\big]\big\rangle
      -\kappa \mu \big\langle \tilde{u}^k,(\e \partial_t)^k
  \big[ \nabla e^{-\e p}\nabla(e^{-\e p+\theta}\nabla \theta)\big]\big\rangle\\
  \leq{}&\frac{l_3}{4}\|\nabla\tilde{u}^k\|^2
          +\frac{\kappa\mu }{l_3}\|e^{-2\e p+\theta}\|^2_{L^{\infty}}\|\nabla \theta\|^2_{s+1,\e}
          +C(\mathcal{M}_1)(1+\mathcal{M}_2)\\
  \leq{}&\frac{l_3}{4}\|\nabla\tilde{u}^k\|^2
          +l_4\|\nabla \theta\|^2_{s+1,\e}
          +C(\mathcal{M}_1)(1+\mathcal{M}_2),
\end{align*}
for some positive constant $l_4$.
Similarly, there exist two positive constants $l_5$ and $l_6$ such that
\begin{align*}
  &Q_4\leq \frac{l_3}{4}\|\nabla\tilde{u}^k\|^2
          +l_5\|\nabla \theta\|^2_{s+1,\e}
          +C(\mathcal{M}_1)(1+\mathcal{M}_2).\\
  &Q_5\leq \frac{l_3}{4}\|\nabla\tilde{u}^k\|^2
          +l_6\|\nabla \theta\|^2_{s+1,\e}
          +C(\mathcal{M}_1)(1+\mathcal{M}_2).
\end{align*}
Finally, the term $Q_6$ and $Q_7$ can be controlled as follows
\begin{align*}
  |Q_6|\leq {}& C(\mathcal{M}_1)(1+\mathcal{M}_2),\\
  |Q_7|\leq {}& C(\mathcal{M}_1)(1+\mathcal{M}_2).
\end{align*}
Putting the estimates of $Q_i$ ($i=1, \dots ,7$) into \eqref{g5} gives
\begin{equation*}
  \big\langle 2\tilde{u}^k,(\e \partial_t)^k g_5\big\rangle
  \leq \frac{3l_3}{4}\|\nabla\tilde{u}^k\|^2+(l_4+l_5+l_6)\|\nabla \theta\|_{s+1,\e}^2+
  C(\mathcal{M}_1)(1+\mathcal{M}_2).
\end{equation*}
For the rest terms in $K_5$, we need only to consider the case $1\leq k\leq s$.
It follows from Corollary \ref{Moser-tx} that
\begin{align*}
\big|\big\langle \tilde{u}^k,-\big [(\e \partial_t)^k,
     e^{-\theta}(\partial_t+ u\cdot \nabla)\big]\tilde{u}\big\rangle\big|
    &  \leq C(\mathcal{M}_1),\\
\big|\big\langle \tilde{u}^k, \big [(\e \partial_t)^k,
      e^{-\e p} \big]\nabla\cdot\Psi(\tilde{u})\big\rangle\big|
 &  \leq C(\mathcal{M}_1)(1+\mathcal{M}_2),\\
\frac{\kappa}{2}\big|\big\langle \tilde{u}^k, \big [(\e \partial_t)^k,
               e^{-\e p} \big]\nabla\nabla \cdot\tilde{u}\big\rangle\big|
         &       \leq C(\mathcal{M}_1)(1+\mathcal{M}_2).
\end{align*}
Then, we obtain that
\begin{equation}\label{k5}
  |K_5|\leq \frac{3l_3}{4}\|\nabla\tilde{u}^k\|^2+(l_4+l_5+l_6)\|\nabla \theta\|_{s+1,\e}^2+
  C(\mathcal{M}_1)(1+\mathcal{M}_2).
\end{equation}

Putting \eqref{k1}-\eqref{k4} and \eqref{k5} into \eqref{dt-p-u} and integrating the resulting inequality yield
\begin{equation}\label{tilda-pkuk}
     \|\tilde{p}^k\|^2+\frac{e^{-\theta}}{2}\|\tilde{u}^k\|^2
       +\frac{l_3}{4}\int_0^T\| \nabla\tilde{u}^k\|^2 \mathrm{d}t
          \leq C_0(M_0)+\sqrt{T}C(\mathcal{M}),
\end{equation}
where we have used \eqref{hat-p-u-theta} to bound $(l_4+l_5+l_6)\int_0^T \|\nabla \theta\|_{s+1,\e}^2\mathrm{d}t$.

Using \eqref{tilda-pkuk}, Lemma \ref{ep-eu-I} and the elementary inequality
\begin{equation*}
  e^{T C(\mathcal{M})}\leq 1+ T \tilde{C}(\mathcal{M}),\;\;\;\forall\,T\in(0,1],
\end{equation*}
we obtain from the expression of the equivalent pressure in \eqref{bianhuan} that
\begin{align*}
  \|(\e \partial_t)^kp\|
  \leq{}& \|\tilde{p}^k\|+\|(\e \partial_t)^k[e^{-\e p}(I_0-I_c)]\|\\
  \lesssim{}& \|\tilde{p}^k\|+\|\e p\|^{s}_{s,\e}\|I_0-I_c\|_{s,\e}\\
  \leq{}& C_0(M_0)+\sqrt[4]{T}C(\mathcal{M})+\big(C_0(M_0)+\sqrt[4]{T}C(\mathcal{M})\big)^{s+1}\\
  \leq{}& C_0(M_0)+\sqrt[4]{T}C(\mathcal{M})+C_0(M_0) e^{(s+1)\sqrt[4]{T}C(\mathcal{M})}\\
  \leq{}& C_0(M_0)+\sqrt[4]{T}C(\mathcal{M}).
\end{align*}
Similarly, we have
\begin{equation*}
  \|(\e \partial_t)^k u\|\lesssim\|\tilde{u}^k\|+
  \|(\e \partial_t)^k(e^{-\epsilon p+\theta}\nabla\theta)\|\leq{} C_0(M_0)+\sqrt[4]{T}C(\mathcal{M}),
\end{equation*}
and
\begin{align*}
  \int_0^T\| \nabla (\e\partial_t)^k u\|^2\mathrm{d}t
  \lesssim{}    &\int_0^T\|\nabla\tilde{u}^k\|^2\mathrm{d}t+
  \int_0^T\|(\e \partial_t)^k\nabla(e^{-\epsilon p+\theta}\nabla\theta)\|^2\mathrm{d}t\\
  \lesssim{} &\int_0^T\|\nabla\tilde{u}^k\|^2\mathrm{d}t+
  \int_0^TC(\mathcal{M})+\|\nabla \theta\|_{s+1,\e}^2\mathrm{d}t\\
  \leq{} & C_0(M_0)+\sqrt{T}C(\mathcal{M}).
\end{align*}
Hence, we conclude that the estimates \eqref{pk-uk} hold.
\end{proof}

Now, we show the estimate of $(\nabla p ,\nabla \cdot u)$.
\begin{lemma}\label{nabla p div u lemma}
Let  $\delta\in (0,2]$ and $s\geq 4$ be an integer. Suppose that $(p,u,\theta,I_0,I_1)$ is a solution to the Cauchy problem \eqref{nondim} and \eqref{initial data} on $[0,T_1]$.
Then for any $\e \in(0,1]$ and $T\in (0,\hat{T}]$, $\hat{T}=\mathrm{min}\{T_1,1\}$, it holds that
\begin{equation}\label{nabla p div u}
   \|(\nabla p,\nabla \cdot u)\|_{s-1,\e,T}
    \leq C_0(M_0)+\big(\e+\sqrt[4]{T}\big)C(\mathcal{M}).
\end{equation}
\end{lemma}

\begin{proof}
We first rewrite \eqref{nondim}$_1$ as
\begin{equation}\label{div u}
    2\nabla \cdot u
    = -\e \partial_t p-\e u \cdot \nabla p
      + e^{-\e p}\big[\kappa \Delta e^{\theta}+\e^2 \Psi(u):\nabla u+\e(I_0-e^{4\theta})
       - (\e^2+\e^{2-\delta})I_1\cdot u\big].
\end{equation}
Applying the operator $(\e \partial_t)^k$, $0\leq k\leq s-1$, to \eqref{div u} and taking the $L^2$-norm of the resulting equation, we have
\begin{align*}
    \|(\e \partial_t)^k\nabla \cdot u\|
    \lesssim{}& \| (\e \partial_t)^{k+1}p\|+ \e\|u\|_{s-1,\e} \| p\|_{s,\e}
        + \|\e p\|_{s-1,\e}^{s-1}\|\theta-\theta_c\|^{s+1}_{s+1,\e}\\
     {}&    +\e\|\e p\|_{s-1,\e}^{s-1}\big[\e\|\nabla u\|^2_{s-1,\e}+\|I_0-I_c\|_{s-1,\e}
          +\|\theta-\theta_c\|^{s-1}_{s-1,\e}  \\
     {}&  +(\e+\e^{1-\delta})\|I_1\|_{s-1,\e}\|u\|_{s-1,\e}\big]\\
     \leq {}& C_0(M_0)+\sqrt[4]{T}C(\mathcal{M})+\e C(\mathcal{M})
          +C_0(M_0)e^{2s\sqrt[4]{T}C(\mathcal{M})}\\
     \leq {}& C_0(M_0)+(\e+\sqrt[4]{T})C(\mathcal{M}),
\end{align*}
where we have used \eqref{pk-uk}, \eqref{deltaI} and Lemma \ref{ep-eu-I}.

Similarly, it follows from \eqref{nondim}$_2$ and \eqref{nondim}$_5$ that
\begin{equation}\label{nabla p}
  \nabla p= -\e e^{-\theta}(\partial_t u+u\cdot \nabla u)
     +\e e^{-\e p}\nabla\cdot \Psi(u)-\frac{e^{-\e p}}{3}(3\e \partial_t I_1+\nabla I_0).
\end{equation}
Applying $(\e \partial_t)^k$, $0\leq k\leq s-1$, to \eqref{nabla p}, we have
\begin{align}
  \|(\e \partial_t)^k\nabla p\|
  \lesssim{}& \|\theta-\theta_c\|^{s-1}_{s-1,\e}\|\e \partial_tu\|_{s-1,\e}
           +\e\|\theta-\theta_c\|^{s-1}_{s-1,\e}\|u\|_{s-1,\e}\|\nabla u\|_{s-1,\e}\notag\\
       {}& +\e\|\e p\|^{s-1}_{s-1,\e}\| \nabla u\|_{s,\e}
           +\|\e p\|^{s-1}_{s-1,\e}\big(\e\|\partial_t I_1\|_{s-1}+\|\nabla I_0\|_{s-1}\big)\notag\\
  \leq{}& C_0(M_0)e^{s\sqrt[4]{T}C(\mathcal{M})}+\e C(\mathcal{M})\notag\\
  \leq{}& C_0(M_0)+(\e+\sqrt[4]{T})C(\mathcal{M})\label{p s-1}.
\end{align}

Furthermore, applying $D^{k,1}$, $0\leq k\leq s-2$, to \eqref{div u} and using \eqref{p s-1}, we obtain that
\begin{equation}\label{div u s-2}
  \|D^{k,1}\nabla \cdot u\| \leq C(M_0)+(\e+\sqrt[4]{T})C(\mathcal{M}).
\end{equation}
It follows from the formula $\nabla \times \nabla=0$ and Lemma \ref{div-curl} that
\begin{equation*}
  \|\nabla \nabla p\|_i
  \leq\|\nabla \cdot \nabla p\|_i, \;\;i=0, \dots ,s-2.
\end{equation*}
Therefore, applying $\nabla \cdot$ to \eqref{nabla p} and using \eqref{div u s-2}, we obtain that
\begin{equation*}
  \|D^{k,1} \nabla p\|\leq\|(\e\partial_t)^k\nabla \cdot \nabla p\|
  \leq C(M_0)+(\e+\sqrt[4]{T})C(\mathcal{M}).
\end{equation*}
Thus, we get the estimates of $\|D^{k,1}(\nabla p,\nabla \cdot u)\|$, $1\leq k\leq s-2$.

Continuing by induction, we finally obtain the desired estimates \eqref{nabla p div u}.
\end{proof}

Next, let us study the estimate of $\nabla \times (e^{-\theta}u)$.
With the help of \eqref{nondim}$_5$, we reformulate the equation \eqref{nondim}$_2$ as
\begin{equation*}
    \partial_t (e^{-\theta}u) + u\cdot \nabla (e^{-\theta}u)
    + \frac{\nabla p}{\epsilon}+e^{-\e p} \frac{\nabla I_0}{3\e}
     = \mu e^{-\e p+\theta} \Delta (e^{-\theta}u) +g_6,
\end{equation*}
where
\begin{align*}
 g_6= {}& -e^{-\e p} \partial_t I_1-e^{-\theta}(\partial_t \theta+u\cdot \nabla \theta)u
         +\mu e^{-\e p}\nabla\cdot\big[ \nabla e^\theta\otimes (e^{-\theta}u)\big]\\
     {}& +\mu e^{-\e p} \nabla(e^{-\theta}u)\nabla e^\theta
         +(\mu+\lambda)e^{-\e p}\nabla \nabla \cdot u.
\end{align*}
Applying the operator $\nabla \times$ to the above equation gives
\begin{equation}\label{curl u}
    (\partial_t  + u\cdot \nabla )\nabla \times(e^{-\theta}u)
     =\mu e^{-\e p+\theta}\Delta [\nabla \times (e^{-\theta}u)]+\nabla \times g_6 +g_7,
\end{equation}
where
\begin{equation*}
  g_7=-\big[\nabla \times,u\cdot \nabla \big](e^{-\theta}u)
        +\mu \big[\nabla \times,  e^{-\e p+\theta}] \Delta (e^{-\theta}u)
        +\frac{1}{3}e^{-\e p}\nabla p \times\nabla I_0.
\end{equation*}
And  we have
\begin{lemma}\label{curl u-lemma}
Let  $\delta\in (0,2]$ and $s\geq 4$ be an integer.
Suppose that $(p,u,\theta,I_0,I_1)$ is a solution to the Cauchy problem \eqref{nondim} and \eqref{initial data} on $[0,T_1]$.
Then for any $\e \in(0,1]$ and $T\in (0,\hat{T}]$, $\hat{T}=\mathrm{min}\{T_1,1\}$, it holds that
\begin{equation}\label{curl-ethetau}
   \|\nabla \times(e^{-\theta}u)\|_{s-1,\e,T}
   +\bigg(\int_0^T\| \nabla\nabla\times
    (e^{-\theta}u)(t)\|^2_{s-1,\e}\mathrm{d}t\bigg)^{\frac{1}{2}}
   \leq C_0(M_0)+\sqrt[4]{T}C(\mathcal{M}).
\end{equation}
\end{lemma}

\begin{proof}
Let $w=\nabla \times(e^{-\theta}u)$.
For any $k$ and $\alpha$ satisfying $0\leq k+|\alpha|\leq s-1$, taking $D^{k,\alpha}$ of \eqref{curl u} yields
\begin{equation*}
  (\partial_t  + u\cdot \nabla )D^{k,\alpha}w
  =\mu e^{-\e p+\theta}\Delta D^{k,\alpha}w + g_8 + g_9,
\end{equation*}
where
\begin{equation*}
   g_8=D^{k,\alpha}\nabla \times \big\{
         \mu e^{-\e p}\nabla\cdot\big[ \nabla e^\theta\otimes (e^{-\theta}u)\big]
         +\mu e^{-\e p} \nabla(e^{-\theta}u)\nabla e^\theta\big\},
\end{equation*}
and
\begin{align*}
   g_9={}& D^{k,\alpha}(\nabla \times g_6+g_7)
           -\big[D^{k,\alpha},u \cdot \nabla\big]w
           +\big[D^{k,\alpha}, \mu e^{-\e p+\theta}]\Delta w
           -g_8\\
      ={}&-D^{k,\alpha}\nabla \times (e^{-\e p}\partial_t I_1)
           -D^{k,\alpha}\nabla \times\big [e^{-\theta}(\partial_t \theta+u\cdot \nabla \theta)u\big]\\
        {}& -(\mu+\lambda) D^{k,\alpha}(\e e^{-\e p} \nabla p \times \nabla \nabla \cdot u)
        -\big[D^{k,\alpha}\nabla \times,u\cdot \nabla \big](e^{-\theta}u)\\
       {}& +\mu\big[D^{k,\alpha}\nabla \times, e^{-\e p+\theta} \big] \Delta (e^{-\theta}u)
        +\frac{1}{3}D^{k,\alpha}(e^{-\e p}\nabla p \times\nabla I_0).
\end{align*}
Multiplying the result with $2D^{k,\alpha}w$ and integrating over $\mathbb{R}^3$, we have
\begin{align}
  \frac{\mathrm{d}}{\mathrm{d}t}\|D^{k,\alpha}w\|^2
   ={}&\int_{\mathbb{R}^3}(\nabla \cdot u)|D^{k,\alpha}w|^2 \mathrm{d}x
   +\big\langle 2D^{k,\alpha}w,\mu e^{-\e p+\theta} \Delta D^{k,\alpha}w \big\rangle \notag\\
   {}& +\big\langle 2D^{k,\alpha}w, g_8\big\rangle
   +\big\langle 2D^{k,\alpha}w, g_9\big\rangle \notag\\
   :={}&\sum\limits_{i=1}^4 R_i.\label{w}
\end{align}
It is straightforward to show that
\begin{equation*}
  |R_1|\leq C(\mathcal{M}_1)\|D^{k,\alpha}w\|^2
\end{equation*}
and
\begin{equation*}
  |R_4|\leq C(\mathcal{M}_1)(1+\mathcal{M}_2).
\end{equation*}
Note that
\begin{equation*}
  \| D^{k,\alpha} w\|\leq \|\theta-\theta_c\|^{s}_{s,\e}\|u\|_{s,\e}\leq C(\mathcal{M}_1)
\end{equation*}
and
\begin{equation*}
  \| D^{k,\alpha}\nabla w\|\leq \|\theta-\theta_c\|^{s+1}_{s+1,\e}
     (\|u\|_{s,\e}+\|\nabla u\|_{s,\e})\leq C(\mathcal{M}_1)(1+\mathcal{M}_2).
\end{equation*}
Thus, for the term $R_2$, we obtain from integration by parts that
\begin{align*}
\big\langle 2D^{k,\alpha}w,\mu e^{-\e p+\theta} \Delta D^{k,\alpha}w \big\rangle
={}& -\big\langle 2D^{k,\alpha}\nabla w,\mu e^{-\e p+\theta}  D^{k,\alpha}\nabla w \big\rangle
    -\big\langle 2D^{k,\alpha} w,\mu   D^{k,\alpha}\nabla w \nabla e^{-\e p+\theta} \big\rangle\\
\lesssim {}&-\|D^{k,\alpha}\nabla w\|^2
         +\|\nabla(\e p,\theta)\|_{L^{\infty}} \| D^{k,\alpha}w\|
              \| D^{k,\alpha}\nabla w\|\\
\leq {}&-\|D^{k,\alpha}\nabla w\|^2+ C(\mathcal{M}_1)(1+\mathcal{M}_2).
\end{align*}
Similarly, the term $R_3$ can be bounded as follows,
\begin{align*}
\big\langle 2D^{k,\alpha}w,g_8\big\rangle
 ={}& \big\langle 2\nabla \times  D^{k,\alpha}w, D^{k,\alpha}\big\{
         \mu e^{-\e p}\nabla\cdot\big[ \nabla e^\theta\otimes (e^{-\theta}u)\big]
         +\mu e^{-\e p} \nabla(e^{-\theta}u)\nabla e^\theta\big\}\\
\lesssim {}& \|\nabla\times D^{k,\alpha}w\|\|\e p\|_{s-1,\e}^{s-1}
         \|\theta-\theta_c\|_{s+1,\e}^{2s}\|u\|_{s,\e}\\
\leq {}&   C(\mathcal{M}_1)(1+\mathcal{M}_2).
\end{align*}
Then, putting the above estimates into \eqref{w} yields
\begin{equation*}
    \frac{\mathrm{d}}{\mathrm{d}t}\|D^{k,\alpha}w\|^2
   +\|D^{k,\alpha}\nabla w\|^2 \leq C(\mathcal{M}_1)(1+\mathcal{M}_2),
\end{equation*}
which implies that \eqref{curl-ethetau} holds.
\end{proof}
Due to the estimates of $\theta$ in \eqref{hat-p-u-theta}, it is easy to show the following corollary.
\begin{corol}\label{corol-curlu}
Under the assumptions in Lemma \ref{curl u-lemma}, it holds that
\begin{equation*}
   \|\nabla \times u\|_{s-1,\e,T}
   +\bigg(\int_0^T\| \nabla\nabla\times u(t)\|^2_{s-1,\e}\mathrm{d}t\bigg)^{\frac{1}{2}}
   \leq C_0(M_0)+\sqrt[4]{T}C(\mathcal{M}).
\end{equation*}
\end{corol}

Up to now, we have obtained the uniform estimate of $\|(p,u)\|_{s,\e}$.
The next task is to close the uniform estimates by showing the estimate of
$\int_0^T\|\nabla (p,u)\|^2_{s,\e}\mathrm{d}t$.
Before that, we need to give the following estimate.
\begin{lemma}
Let  $\delta\in (0,2]$ and $s\geq 4$ be an integer.
Suppose that $(p,u,\theta,I_0,I_1)$ is a solution to the Cauchy problem \eqref{nondim} and \eqref{initial data} on $[0,T_1]$.
Then for any $\e \in(0,1]$ and $T\in (0,\hat{T}]$, $\hat{T}=\mathrm{min}\{T_1,1\}$, it holds that
\begin{equation}\label{p s+1}
  \int_0^T\|(\e \partial_t)^{s+1} p\|^2\mathrm{d}t
    \leq C_0(M_0)+\sqrt{T}C(\mathcal{M}).
\end{equation}
\end{lemma}
\begin{proof}
Applying the operator $\e(\e\partial_t)^s$ to \eqref{nondim}$_1$ yields
\begin{align}
    (\e \partial_t)^{s+1} p
    = (\e \partial_t)^{s}\big\{
    &  -u \cdot \nabla (\e p)-2\nabla \cdot u
       +\kappa e^{-\e p+\theta}(\nabla\theta\cdot \nabla \theta+\Delta \theta )\notag\\
    &  + e^{-\e p}\big[ \Psi(\e u):\nabla (\e u)+\e(I_0-e^{4\theta})
       - (\e^2+\e^{2-\delta})  I_1\cdot u\big]
         \big\}.
\end{align}
And then, we have
\begin{align*}
  \|(\e \partial_t)^{s+1} p\|
   \lesssim {} & \| u\|_{s,\e}\|\e p\|_{s+1,\e}
                 +\|\nabla(\e \partial_t)^{s} u\|
                   +\|(\e p, \theta-\theta_c)\|^{s}_{s,\e}\|\nabla \theta\|^2_{s,\e}\\
               &   +\|(\e p, \theta-\theta_c)\|^{s}_{s,\e}\|\nabla \theta\|_{s,\e}
                   +\|e^{-\e p+\theta}\|_{L^{\infty}}\|\nabla \theta\|_{s+1,\e}\\
              & +\|\e p\|^s_{s,\e}\big(\|\nabla (\e u)\|^2_{s,\e}
                      +\e\|I_0-I_c\|_{s,\e}+\e\|\theta-\theta_c\|^s_{s,\e}
                       +\|I_1\|_{s,\e}\|u\|_{s,\e}\big)  \\
    \lesssim {} &  \|\nabla(\e \partial_t)^{s} u\|+\|\nabla\theta\|_{s+1,\e}
                     +C(\mathcal{M}_1),
\end{align*}
which implies that \eqref{p s+1} holds by using \eqref{hat-p-u-theta} and \eqref{pk-uk}.
\end{proof}

\begin{lemma}\label{int nabla p u lemma}
Let  $\delta\in (0,2]$ and $s\geq 4$ be an integer.
Suppose that $(p,u,\theta,I_0,I_1)$ is a solution to the Cauchy problem \eqref{nondim} and \eqref{initial data} on $[0,T_1]$.
Then for any $\e \in(0,1]$ and $T\in (0,\hat{T}]$, $\hat{T}=\mathrm{min}\{T_1,1\}$, it holds that
\begin{equation}
\int_0^T\|\nabla (p,u)\|^2_{s,\e}\mathrm{d}t
    \leq C_0(M_0)+\sqrt{T}C(\mathcal{M}).
\end{equation}
\end{lemma}

\begin{proof}
According to the fact $\|(p ,u)\|_{s,\e}\leq \mathcal{M}_1$ and Corollary \ref{corol-curlu}, it is sufficient to show that
\begin{equation}\label{nablap-divu-s+1}
\int_0^T\|D^{s-i,i}(\nabla p,\nabla \cdot u)\|^2\mathrm{d}t
    \leq C_0(M_0)+\sqrt{T}C(\mathcal{M}),\;\; i=0,1, \dots , s.
\end{equation}

First of all, we consider the case $i=0$.
Taking $(\e \partial_t)^s$ of \eqref{nabla p} yields
\begin{equation*}
  (\e \partial_t)^s\nabla p
  = -(\e \partial_t)^s( e^{-\theta}\e\partial_t u)
    +(\e \partial_t)^s \bigg[-\e e^{-\theta}u\cdot \nabla u +\e e^{-\e p}\nabla\cdot\Psi(u)
      -e^{-\e p}\bigg(\e \partial_t I_1+\frac{\nabla I_0}{3}\bigg)\bigg],
\end{equation*}
and we then obtain from Cauchy's inequality and integration by parts that
\begin{align*}
  \|(\e \partial_t)^s\nabla p\|^2
\leq{} & -\big\langle(\e \partial_t)^s\nabla p,
          (\e \partial_t)^s( e^{-\theta}\e\partial_t u)\big\rangle
              +\frac{1}{2}\|(\e \partial_t)^s\nabla p\|^2\\
   {}& +\frac{1}{2}\|(\e \partial_t)^s[- \e e^{-\theta}u\cdot \nabla u
          +  \e e^{-\e p}\nabla\cdot \Psi( u)
          -e^{-\e p}(\e \partial_t I_1+3^{-1}\nabla I_0)]\|^2\\
  \leq {} & -\e\frac{\mathrm{d}}{\mathrm{d}t}
             \big\langle(\e \partial_t)^s\nabla p,
             (\e \partial_t)^{s-1}( e^{-\theta}\e\partial_t u)\big\rangle
            +\big\langle(\e \partial_t)^{s+1}\nabla p,
              (\e \partial_t)^{s-1}( e^{-\theta}\e\partial_t u)\big\rangle\\
       {}& +\frac{1}{2}\|(\e \partial_t)^s\nabla p\|^2
           +C(\mathcal{M}_1)+l_7 \|\nabla(\e u)\|^2_{s+1,\e}\\
  \leq {} & -\e\frac{\mathrm{d}}{\mathrm{d}t}
             \big\langle(\e \partial_t)^s\nabla p,
             (\e \partial_t)^{s-1}( e^{-\theta}\e\partial_t u)\big\rangle
            -\big\langle(\e \partial_t)^{s+1} p,
              (\e \partial_t)^{s-1}\nabla\cdot( e^{-\theta}\e\partial_t u)\big\rangle\\
       {}& +\frac{1}{2}\|(\e \partial_t)^s\nabla p\|^2
           +C(\mathcal{M}_1)+l_7\|\nabla(\e u)\|^2_{s+1,\e}\\
  \leq {} & -\e\frac{\mathrm{d}}{\mathrm{d}t}
             \big\langle(\e \partial_t)^s\nabla p,
             (\e \partial_t)^{s-1}( e^{-\theta}\e\partial_t u)\big\rangle
            +\frac{1}{2}\|(\e \partial_t)^{s+1} p\|^2\\
       {}& +\frac{1}{2}\|(\e\partial_t)^{s}\nabla \cdot u\|^2
           +\frac{1}{2}\|(\e \partial_t)^s\nabla p\|^2
           +C(\mathcal{M}_1)+l_7\|\nabla(\e u)\|^2_{s+1,\e},
\end{align*}
where $l_7$ is a positive constant.
Integrating the above inequality on $[0,T]$ and using \eqref{hat-p-u-theta}, \eqref{pk-uk} and \eqref{p s+1}, we obtain that
\begin{align}
  \frac{1}{2}\int_0^T\|(\e \partial_t)^s\nabla p\|^2 \mathrm{d}t
    \leq {} & -\big\langle(\e \partial_t)^s\nabla (\e p),
             (\e \partial_t)^{s-1}( e^{-\theta}\e\partial_t u)\big\rangle\big|^T_0
             +\frac{1}{2}\int_0^T\|(\e \partial_t)^{s+1} p\|^2\mathrm{d}t \notag\\
        {}&  +\frac{1}{2}\int_0^T\|(\e\partial_t)^{s}\nabla\cdot u\|^2\mathrm{d}t
             +TC(\mathcal{M}_1)
             +l_7\int_0^T\|\nabla(\e u)\|^2_{s+1,\e}\mathrm{d}t  \notag\\
     \leq {} & \sup\limits_{t\in[0,T]}
              \big\{\|\e p\|_{s+1,\e}\|\theta-\theta_c\|^{s-1}_{s-1,\e}\|u\|_{s,\e}\big\}
             +C_0(M_0)+\sqrt{T}C(\mathcal{M}) \notag\\
    \leq {} & C_0(M_0)+\sqrt{T}C(\mathcal{M}).\label{int p}
\end{align}
Hence, we arrive at
\begin{equation*}
\int_0^T\|(\e \partial_t)^s(\nabla p,\nabla \cdot u)\|^2\mathrm{d}t
    \leq C_0(M_0)+\sqrt{T}C(\mathcal{M}).
\end{equation*}

Next, for the case $i=1$, applying the operator $D^{s-1,1}$  to \eqref{div u}, we have
\begin{equation*}
    \|D^{s-1,1}\nabla \cdot u\|
    \lesssim \|D^{s,1} p\|+\|\nabla\theta\|_{s+1,\e}+C(\mathcal{M}_1).
\end{equation*}
Using \eqref{hat-p-u-theta} and \eqref{int p}, we obtain that
\begin{align}
\int_0^T\|D^{s-1,1}\nabla \cdot u\|^2\mathrm{d}t
\lesssim {}& \int_0^T \|D^{s,1} p\|^2+\|\nabla\theta\|^2_{s+1,\e}+C(\mathcal{M}_1)\mathrm{d}t\notag\\
    \leq {} &C_0(M_0)+\sqrt{T}C(\mathcal{M}).\label{Du,s-1}
\end{align}
Similarly, applying the operator $D^{s-1,1}$  to \eqref{nabla p}, we obtain that
\begin{equation*}
  \|D^{s-1,1}\nabla p\|\lesssim  \|D^{s,1} u\|
        +\|\nabla(\e u)\|_{s+1,\e}+C(\mathcal{M}_1).
\end{equation*}
Then, it follows that
\begin{align}
\int_0^T\|D^{s-1,1}\nabla p\|^2\mathrm{d}t
\lesssim {}& \int_0^T \|D^{s,1} u\|^2
  +\|\nabla(\e u)\|^2_{s+1,\e}+C(\mathcal{M}_1) \mathrm{d}t\notag\\
    \leq {} &C_0(M_0)+\sqrt{T}C(\mathcal{M}).\label{Dp,s-1}
\end{align}
Combining \eqref{Du,s-1} and \eqref{Dp,s-1}, we arrive at
\begin{equation*}
\int_0^T\|D^{s-1,1}(\nabla p,\nabla \cdot u)\|^2\mathrm{d}t
    \leq C_0(M_0)+\sqrt{T}C(\mathcal{M}).
\end{equation*}

By applying the operator $D^{s-i,i}$, $i=2, \dots ,s$, to \eqref{div u} and \eqref{nabla p} in turn, we finally obtain \eqref{nablap-divu-s+1} from induction, and hence complete the proof.
\end{proof}

\begin{proof}[Proof of Proposition\ref{priori}]
Proposition \ref{priori} follows directly from Lemmas \ref{ep-eu-I}, \ref{pt-ut-lemma},
\ref{nabla p div u lemma} and \ref{int nabla p u lemma}, and Corollary \ref{corol-curlu}.
\end{proof}

\section{Convergence for the case $\delta \in (0,2]$}\label{S.4}
In this section, we shall prove Theorem \ref{convergence thm} by modifying the arguments developed by M\'{e}tivier and Schochet \cite{MS}, which mainly contains the local energy decay of the acoustic wave equations and the method of compactness argument.

First of all, we give the following lemma which describes the convergence of slow components of the solutions, and can be obtained from the uniform estimates
\eqref{uniform estimate}.
\begin{lemma}\label{converge 1}
Suppose that the assumptions stated in Theorem \ref{uniform existence theorem} hold.
Then there exist a quintuple $(\bar{p},\bar{u},\bar{\theta},\bar{I}_0,\bar{I}_1)$ satisfying
$(\bar{p},\bar{u})\in L^{\infty}(0,T_0;H^s(\mathbb{R}^3))$ and
$(\bar{\theta}-\theta_c,\bar{I}_0-I_c,\bar{I}_1)\in L^{\infty}(0,T_0;H^{s+1}(\mathbb{R}^3))$ such that, after extracting a subsequence,
\begin{equation}\label{p,u infty}
  (p^\e,u^\e)\rightarrow (\bar{p},\bar{u})\;\,\;
  \mathrm{weakly-*}\;\mathrm{in}\;\; L^{\infty}(0,T_0;H^s(\mathbb{R}^3)),
\end{equation}
\begin{equation}\label{theta I0 I1 infty}
  (\theta^\e-\theta_c,I^\e_0-I_c,I^\e_1)
  \rightarrow(\bar{\theta}-\theta_c,\bar{I}_0-I_c,\bar{I}_1)\;\,\;
  \mathrm{weakly-*}\;\mathrm{in}\;\, L^{\infty}(0,T_0;H^{s+1}(\mathbb{R}^3)).
\end{equation}
Moreover, after further extracting a subsequence, we have
\begin{align}
  (I^\e_0-I_c,I^\e_1)\rightarrow(\bar{I}_0-I_c,\bar{I}_1)\;\,
  & \mathrm{strongly}\;\mathrm{in}\;\, C([0,T_0];H^{s'+1}_{\mathrm{loc}}(\mathbb{R}^3)),
  \label{I limit}\\
  \theta^\e-\theta_c\rightarrow\bar{\theta}-\theta_c\;\,
  &\mathrm{strongly}\;\mathrm{in}\;\, C([0,T_0];H^{s'+1}_{\mathrm{loc}}(\mathbb{R}^3)),
  \label{theta limit}\\
   \nabla\times (e^{-\theta^{\e}}u^\e)\rightarrow
   \nabla\times (e^{-\bar{\theta}}\bar{u})\;\,
   &\mathrm{strongly}\;\mathrm{in}\;\, C([0,T_0];H^{s'-1}_{\mathrm{loc}}(\mathbb{R}^3)),
   \label{curl u limit}
\end{align}
for all $s'<s$.
\end{lemma}

\begin{proof}
The uniform estimates \eqref{uniform estimate} gives us
\begin{equation}\label{bounded}
   \|(p^\e,u^\e)\|_{s,\e,T_{0}}+\| \theta^\e-\theta_c\|_{s+1,\e,T_0}
   +\interleave (I^\e_0-I_c,I^\e_1)\interleave_{s+1,\e,T_0}<M_1,
\end{equation}
which implies that, after extracting subsequences, \eqref{p,u infty} and \eqref{theta I0 I1 infty} hold.

In addition, according to the definition of the weighted norm $\interleave (I^\e_0-I_c,I^\e_1)\interleave_{s+1,\e,T_0}$, it follows that
\begin{equation*}
   (\partial_t I^\e_0,\partial_tI^\e_1)\in  L^{\infty}(0,T_0;H^s(\mathbb{R}^3)),
\end{equation*}
and then we get \eqref{I limit} by using Aubin-Lions Lemma (see \cite{Simon}).

From the equation of $\theta^\e$ in \eqref{nondim}$_3$, we find that
\begin{equation}\label{theta-t infty}
  \partial_t \theta^\e \in C([0,T_0];H^{s-1}(\mathbb{R}^3)).
\end{equation}
Combining \eqref{bounded} and \eqref{theta-t infty}, we obtain \eqref{theta limit}.
Similarly, from \eqref{curl u}, we have
\begin{equation*}
  \partial_t \nabla\times (e^{-\theta^{\e}}u^\e)\in C([0,T_0];H^{s-3}(\mathbb{R}^3)),
\end{equation*}
and then \eqref{curl u limit} follows.
\end{proof}

Next, we show the convergence of fast components of the solution by using the local energy decay of  acoustic wave equations.
\begin{lemma}\label{B}
Suppose that the assumptions stated in Theorem \ref{convergence thm} hold.
Then, for all $s'<s$, we have
\begin{align}
  p^\e+\frac{e^{-\e p^\e}(I^\e_0-I_c)}{3}
  & \rightarrow 0 \;\,\mathrm{strongly}\;\mathrm{in}\;\, L^2([0,T_0];H_{\mathrm{loc}}^{s'}(\mathbb{R}^3)),\label{p converge}\\
 \nabla \cdot (2u^{\epsilon}-\kappa e^{-\e p^\e+\theta^\e}\nabla \theta^\e)
 & \rightarrow 0 \;\,\mathrm{strongly}\;\mathrm{in}\;\, L^2([0,T_0];H_{\mathrm{loc}}^{s'-1}(\mathbb{R}^3)).
 \label{div-thetau converge}
\end{align}
\end{lemma}

\begin{proof}

Recall the definitions of the equivalent pressure $\tilde{p}$ and
the equivalent velocity $\tilde{u}$ in \eqref{bianhuan}, and the equations \eqref{pt-ut} of them.
Applying $\e^2\partial_t$ to \eqref{pt-ut}$_1$ yields
\begin{equation}\label{ptt}
     \e^2\partial^2_t \tilde{p}
     +\e^2\partial_t (u^\e\cdot \nabla \tilde{p})
     +\e\partial_t\nabla \cdot \tilde{u}
     =\e^2\partial_t g_4.
\end{equation}
Multiplying \eqref{pt-ut}$_2$ by $2 e^{\theta^\e}$ and then applying $\e \nabla \cdot$ to the resulting equation yields
\begin{equation}\label{div-ut}
   \e \partial_t \nabla \cdot \tilde{u}
       + \e \nabla \cdot (u^\e\cdot \nabla \tilde{u})
       +  \nabla \cdot (2e^{\theta^\e}\nabla \tilde{p})
   =\e \nabla \cdot \Big[e^{-\e p^\e+\theta^\e} \Big(\nabla\cdot \Psi(\tilde{u})+\frac{\kappa}{2}\nabla \nabla \cdot\tilde{u}\Big)+ 2 e^{\theta^\e}g_5\Big].
\end{equation}
Subtracting \eqref{div-ut} from \eqref{ptt}, we obtain that
\begin{equation*}
  \e^2\partial^2_t \tilde{p}- \nabla \cdot (2e^{\theta^\e}\nabla \tilde{p})=F,
\end{equation*}
where
\begin{equation*}
  F=-\e^2\partial_t (u^\e\cdot \nabla \tilde{p})
       + \e \nabla \cdot (u^\e\cdot \nabla \tilde{u})
        +\e^2\partial_t g_4
        -\e \nabla \cdot \Big[e^{-\e p^\e+\theta^\e} \Big(\nabla\cdot \Psi(\tilde{u})+\frac{\kappa}{2}\nabla \nabla \cdot\tilde{u}\Big)+ 2 e^{\theta^\e}g_5\Big].
\end{equation*}
By virtue of the uniform boundedness of $(p^\e,u^\e,\theta^\e,I^\e_0,I^\e_1)$, we have
\begin{equation*}
   F \rightarrow 0 \;\,\mathrm{strongly}\;\mathrm{in}\;\, L^2([0,T_0];L^{2}(\mathbb{R}^3)).
\end{equation*}
According to \eqref{decay theta}, the rapid decay condition of $\theta^\e_0$ at infinity, and the strong convergence of $\theta^\e$, it is easy to show the coefficient $2e^{\theta^\e}$ in $\nabla \cdot (2e^{\theta^\e}\nabla p^\e)$ satisfies the requirement of $b^\e$ in Lemma \ref{dispersive}.
Therefore, from Lemma \ref{dispersive}, we get
\begin{equation*}
  \tilde{p}\rightarrow 0 \;\,\mathrm{strongly}\;\mathrm{in}\;\,
   L^2([0,T_0];L^2_{\mathrm{loc}}(\mathbb{R}^3)).
\end{equation*}
Since we had the boundedness of $ (p^\e, I^\e_0-I_c)$ in $L^{\infty}([0,T_0];H^s(\mathbb{R}^3))$,
an interpolation argument gives \eqref{p converge}.

Similarly, setting $\phi^\e=\nabla \cdot \tilde{u}$,
the acoustic wave equation for $\phi^\e$ is given by
\begin{equation*}
  \e^2\partial^2_t \phi^\epsilon- \nabla \cdot (2e^{\theta^\e}\nabla \phi^\e)=G,
\end{equation*}
where
\begin{align*}
  G={}&-\e^2\partial_t \nabla \cdot(u^\e \cdot \nabla \tilde{u})
    -\e \nabla\cdot\big(2e^{\theta^\e}\partial_t \theta^\e\nabla \tilde{p}\big)
    -\e \nabla\cdot\big[2e^{\theta^\e}\nabla(u^\e \cdot \nabla\tilde{p})\big]\\
    &+\e^2\partial_t\nabla \cdot \Big[e^{-\e p^\e+\theta^\e} \Big(\nabla\cdot \Psi(\tilde{u})+\frac{\kappa}{2}\nabla \nabla \cdot\tilde{u}\Big)+ 2 e^{\theta^\e}g_5\Big]
    -\e \nabla\cdot\big(2e^{\theta^\e}\nabla g_4\big),
\end{align*}
and it is easy to see that
$G\rightarrow 0$ in $L^2([0,T_0];L^{2}(\mathbb{R}^3))$.
Therefore, we obtain that
\begin{equation*}
  \nabla \cdot (2u^{\epsilon}-\kappa e^{-\e p^\e+\theta^\e}\nabla \theta^\e) \rightarrow 0 \;\,\mathrm{strongly}\;\mathrm{in}\;\, L^2([0,T];L_{\mathrm{loc}}^2(\mathbb{R}^3)),
\end{equation*}
and \eqref{div-thetau converge} follows.
\end{proof}

Now, we are in a position to show Theorem \ref{convergence thm}.

\begin{proof}[Proof of Theorem \ref{convergence thm}]
Due to  \eqref{theta limit}, \eqref{curl u limit} and \eqref{div-thetau converge}, we obtain that
\begin{equation*}
  \nabla \cdot u^\e \rightarrow \nabla \cdot \bar{u} \;\,\mathrm{strongly}\;\mathrm{in}\;\, L^2([0,T_0];H_{\mathrm{loc}}^{s'-1}(\mathbb{R}^3)),
\end{equation*}
\begin{equation*}
  \nabla \times u^\e \rightarrow \nabla \times \bar{u} \;\,\mathrm{strongly}\;\mathrm{in}\; \, L^2([0,T_0];H_{\mathrm{loc}}^{s'-1}(\mathbb{R}^3)),
\end{equation*}
and then
\begin{equation*}
  u^\e \rightarrow \bar{u} \;\,\mathrm{strongly}\;\mathrm{in}\; \, L^2([0,T_0];H_{\mathrm{loc}}^{s'}(\mathbb{R}^3)).
\end{equation*}
Furthermore, it follows from
 \eqref{p converge} and the convergence of $I_0$ in \eqref{I limit} that
\begin{equation}\label{I1 converge}
  p^\e \rightarrow -\frac{(\bar{I}_0-I_c)}{3} \;\,\mathrm{strongly}\;\mathrm{in}\;\, L^2([0,T_0];H_{\mathrm{loc}}^{s'}(\mathbb{R}^3)).
\end{equation}
Thus, we obtain the desired convergence of $(p^\e,u^\e,\theta^\e-\theta_c,I_0^\e-I_c,I_1^\e)$.

Let us derive the limit system satisfied by $(\bar{p},\bar{u},\bar{\theta},\bar{I}_0,\bar{I}_1)$.

We first consider the case $\delta=2$.
Multiplying \eqref{nondim}$_5$ by $\e^2$ and then taking $\e\rightarrow 0$ yields
\begin{equation*}
   \bar{I}_1=0.
\end{equation*}
Applying the operator $\e \nabla \cdot$ to \eqref{nondim}$_5$, we have
\begin{equation}\label{e1}
  \frac{\nabla\cdot I^\e_1}{\e}
  =-\e \partial_t \nabla\cdot I^\e_1-\frac{\Delta I^\e_0}{3}-\e\nabla\cdot I^\e_1.
\end{equation}
By substituting \eqref{e1} into \eqref{nondim}$_4$, we get
\begin{equation*}
  \partial_t I^\epsilon_0
    -\frac{\Delta I^\e_0}{3}=e^{4\theta^\e}-I^\epsilon_0 +\e \partial_t \nabla\cdot I^\e_1
    +\e\nabla\cdot I^\e_1,
\end{equation*}
and then take $\e\rightarrow 0$ in the above equation to obtain the diffusion equation
\begin{equation}\label{I0 limit}
   \partial_t \bar{I}_0 - \frac{\Delta \bar{I}_0}{3}=e^{4\bar{\theta}}-\bar{I}_0.
\end{equation}

On the other hand, since \eqref{I1 converge} and the fact that $\bar{I}_1=0$, we multiply \eqref{nondim}$_1$ by $\e$ and then pass to the limit in the resulting equation to obtain that
\begin{equation}\label{limit1}
  \nabla \cdot(2\bar{u}- \kappa e^{\bar{\theta}}\nabla\bar{\theta} )=0,
\end{equation}
in the sense of distributions.
Passing to the limit in \eqref{nondim}$_3$, we find, in the sense of distributions, that
\begin{equation}\label{limit2}
 \partial_t \bar{\theta} + \bar{u}\cdot \nabla \bar{\theta}+ \nabla \cdot \bar{u}
    =   \Delta e^{\bar{\theta}}.
\end{equation}
Moreover, adding \eqref{nondim}$_5$ to \eqref{nondim}$_2$, and then applying the operator $\nabla \times$ to the resulting equation, one has
\begin{equation*}
  \nabla \times \big[ e^{-\theta^\e}(\partial_t u^\epsilon + u^\epsilon\cdot \nabla u^\epsilon)
    -e^{-\e p^\e}\nabla\cdot \Psi(u^\e)+e^{-\e p^\e}\partial_t I^\e_1\big]
    -e^{-\e p^\e}\frac{\nabla p^\e\times\nabla I^\e_0}{3}=0.
\end{equation*}
Passing to the limit in the above equations and using \eqref{I1 converge}, we find that
\begin{equation*}
  \nabla \times \big[e^{-\bar{\theta}} (\partial_t \bar{u} + \bar{u}\cdot \nabla \bar{u})
  -\nabla\cdot \Psi(\bar{u})\big]=0
\end{equation*}
holds in the sense of distributions,
which means that
\begin{equation}\label{limit3}
  e^{-\bar{\theta}} (\partial_t \bar{u} + \bar{u}\cdot \nabla \bar{u})+\nabla \pi_1
  =\nabla\cdot \Psi(\bar{u}),
\end{equation}
for some function $\pi_1$.

In addition, according to  \eqref{curl u limit} and \eqref{div-thetau converge}, we deduce, by using the same argument as that in the proof of Theorem 1.5 in \cite{MS}, that the initial data of $(\bar{u},\bar{\theta},\bar{I}_{0})$ are given by
\begin{equation}\label{limit data}
  (\bar{u},\bar{\theta},\bar{I}_{0})|_{t=0}=(\bar{w}_0,\bar{\theta}_0,\bar{I}_{00}),
\end{equation}
 where $\bar{w}_0$ is determined by
\begin{equation*}
  \nabla \cdot(2\bar{w}_0- \kappa e^{\bar{\theta}_0}\nabla \bar{\theta}_0)=0\;\;\mathrm{and}\;\;
  \nabla\times (e^{-\bar{\theta}_0}\bar{w}_0)=\nabla\times (e^{-\bar{\theta}_0}\bar{u}_0).
\end{equation*}
Finally, a standard iterative method shows that the limit system \eqref{I0 limit}-\eqref{limit3} with the initial data \eqref{limit data} has a unique solution $(\bar{u},\bar{\theta},\bar{I}_{0})\in C([0,T_0];H^s(\mathbb{R}^3))$, which implies that the above convergence holds for the full sequence of $(p^\e,u^\e, \theta^\e, I_0^\e,I_1^\e)$.

Next, we focus on the case  $\delta\in (1,2)$.
Multiplying \eqref{nondim}$_5$ by $\e^\delta$ and then taking $\e\rightarrow 0$ yields
\begin{equation*}
   \bar{I}_1=0,
\end{equation*}
and then applying the operator $\e^{\delta-1} \nabla \cdot$ to \eqref{nondim}$_5$ yields
\begin{equation*}
  \frac{\nabla\cdot I^\e_1}{\e}
  =-\e^{\delta-1}\partial_t \nabla\cdot I^\e_1
    -\frac{\Delta I^\e_0}{3\e^{2-\delta}}-\e^{\delta-1}\nabla\cdot I^\e_1.
\end{equation*}
By substituting the above equation into \eqref{nondim}$_4$ and taking $\e\rightarrow 0$, we deduce that
\begin{equation*}
  \Delta \bar{I}_0=0.
\end{equation*}
The derivation of the limit equations satisfied by $(\bar{u},\bar{\theta})$ is as same as in the case of $\delta=2$ and then is omitted.

Finally, we can follow a similar procedure as above to deal with the cases   $\delta=1$ and $\delta\in (0,1)$.
Hence we complete the proof of Theorem \ref{convergence thm}.
\end{proof}

\section{Proof of Theorem \ref{Ls 0 thm}}\label{S.5}
The purpose of this section is to prove Theorem \ref{Ls 0 thm}, which stated the uniform estimates and convergence of the solutions to the system \eqref{nondim} with $\delta=0$.
The idea of proving outline of Theorem \ref{Ls 0 thm} is essentially similar to those of Theorems \ref{uniform existence theorem} and \ref{convergence thm}.
Here, we only give some explanations and point out how to modify them to be applied to  Theorem
\ref{Ls 0 thm}.

Based on the general initial data condition \eqref{initial condition Ls0} and the system \eqref{nondim} with $\delta=0$, we introduce the quantity
\begin{equation*}
  \mathfrak{M}(T)=\sup\limits_{t\in[0,T]}\mathfrak{M}_1(t)
  +\bigg(\int_0^T \mathfrak{M}^2_2(t)\mathrm{d}t \bigg)^{\frac{1}{2}},
\end{equation*}
where
\begin{equation*}
  \mathfrak{M}_1(t)=\|(p^\e,u^\e)(t)\|_{s,\e}
     +\|(\e p^\e,\e u^\e, \theta^\e-\theta_c)(t)\|_{s+1,\e}
     +\|(I^\e_0-I_c,I^\e_1)(t)\|_{s+1,\e}
\end{equation*}
and
\begin{equation*}
  \mathfrak{M}_2(t)=\|\nabla(p^\e,u^\e)(t)\|_{s,\e}
      +\|\nabla(\e u^\e,\theta^\e)(t)\|_{s+1,\e}.
\end{equation*}
Similar to Proposition \ref{priori}, to establish the uniform existence of the solutions, it is sufficient to show that, for any $\e \in(0,\e'_0]$ and $T\in (0,\hat{T}']$, it holds that
\begin{equation*}
  \mathfrak{M}(T)\leq C_0(M_0')+(\sqrt[4]{T}+\e)C(\mathfrak{M}(T)),
\end{equation*}
which can be obtained from the same procedure as that in Section \ref{S.3}.
Then we get the uniform estimates \eqref{Ls-uniform estimate} and the uniform existence of solutions.
We remark that, in this process, there is no need to introduce the equivalent pressure $\tilde{p}$, but only the equivalent velocity $\tilde{u}$.

Since the convergence of $(I^\e_0,I^\e_1)$ is not implied by \eqref{Ls-uniform estimate} directly,
we are mainly concerned with how to obtain the compactness of $(I^\e_0,I^\e_1)$.
Applying $\nabla \times$ to \eqref{nondim}$_5$, we get
\begin{equation*}
  \partial_t\nabla\times I^\e_1 =-2 I^\e_1 \in C([0,T'_0];H^{s-1}(\mathbb{R}^3)).
\end{equation*}
It follows from Aubin-Lions Lemma that
\begin{equation*}
  \nabla\times I_1 \rightarrow \nabla\times \bar{I}_1\;\,
   \mathrm{strongly}\;\mathrm{in}\;\, C([0,T'_0];H^{s'}_{\mathrm{loc}}(\mathbb{R}^3)),
\end{equation*}
for all $s'<s$.
Using the same argument as that in the proof of Lemma \ref{B}, we obtain that
\begin{equation*}
  I_0^\e-I_c \rightarrow 0 \;\,\mathrm{strongly}\;\mathrm{in}\; \, L^2([0,T'_0];H_{\mathrm{loc}}^{s'+1}(\mathbb{R}^3))
\end{equation*}
and
\begin{equation*}
  \nabla \cdot I^{\epsilon}_1\rightarrow 0 \;\,\mathrm{strongly}\;\mathrm{in}\;\, L^2([0,T'_0];H_{\mathrm{loc}}^{s'}(\mathbb{R}^3)).
\end{equation*}
Hence, we get
\begin{equation*}
  (I_0^\e-I_c,I_1^\e) \rightarrow (0,\bar{I}_1) \;\,\mathrm{strongly}\;\mathrm{in}\; \, L^2([0,T'_0];H_{\mathrm{loc}}^{s'+1}(\mathbb{R}^3)).
\end{equation*}
In addition, the convergence of $(p^\e,u^\e,\theta^\e)$ to $(0,\bar{u},\bar{\theta})$ can be obtained by the similar method to that employed in the previous section.
Thus, we achieve our aim to get the strong convergence of the solutions for the case $\delta=0$.

The process of deriving the limit equations \eqref{limit eq-3} satisfied by $(\bar{u},\bar{\theta},\bar{I}_1)$ is analogous to that in Theorem \ref{convergence thm} and we will not repeat it here.
This completes the proof of Theorem \ref{Ls 0 thm}. $\hfill \square$

\begin{appendices}
\section{Appendix}
The purpose of this appendix is to derive the dimensionless equations \eqref{NSFP1-Ma} for the fluids obeying the perfect gas relations \eqref{gas}.
We introduce the following dimensionless quantities:
\begin{gather*}
    x_*=\frac{x}{L_{\infty}},\;\;
   t_*=\frac{t}{T_{\infty}},\;\;
   u_*=\frac{u}{u_{\infty}},\;\;
   \rho_*=\frac{\rho}{\rho_{\infty}},\;\;
   \Theta_*=\frac{\Theta}{\Theta_{\infty}},\;\;
   P_*=\frac{P}{P_{\infty}}, \;\;
   e_*=\frac{e}{e_{\infty}},\\
    I_{0,*}=\frac{I_0}{I_{\infty}},\;\;
   I_{1,*}=\frac{I_1}{I_{\infty}},\;\;
    \mu_*=\frac{\mu}{\mu_{\infty}},\;\;
    \lambda_*=\frac{\lambda}{\mu_{\infty}},\;\;
   \kappa_*=\frac{\kappa}{\kappa_{\infty}},\;\;
   \sigma_{a,*}=\frac{\sigma_{a}}{\sigma_{a,\infty}},\;\;
   \sigma_{s,*}=\frac{\sigma_{s}}{\sigma_{s,\infty}}.
\end{gather*}
Here the symbol with the subscript ``$\infty$" denotes the corresponding characteristic value.
Putting the above scalings into \eqref{NSFP1} and using the perfect gas relations \eqref{gas} and the following compatibility relations:
\begin{equation*}
  e_{\infty}=c_V\Theta_{\infty},\;\;\; P_{\infty}=R\rho_{\infty} \Theta_{\infty},\;\;
  I_{\infty} =c a_r \Theta^4_{\infty},
\end{equation*}
we obtain the dimensionless form of \eqref{NSFP1} as follows
\begin{equation}\label{A.1}
\left\{
\begin{aligned}
  & \mathrm{St}\partial _t \rho  +  u\cdot \nabla \rho+ \rho\nabla \cdot  u = 0, \\
  & \mathrm{St}\rho \partial _t u + \rho u\cdot \nabla u + \frac{1}{\gamma \mathrm{Ma^2}}\nabla P
      =\frac{1}{\mathrm{Re}}\nabla\cdot \Psi(u)
        +\frac{\mathcal{P}\mathcal{L}}{\gamma(\gamma-1) \mathrm{Ma^2}}
        (\sigma_a+\mathcal{L}_s\sigma_s)I_1,\\
  & \mathrm{St}\rho \partial _t \Theta + \rho u\cdot \nabla \Theta+(\gamma-1) P\nabla \cdot u\\
  &  \quad =\frac{\gamma}{\mathrm{RePr}}\kappa\Delta \theta
      +\frac{\gamma(\gamma-1)\mathrm{Ma}^2}{\mathrm{Re}}\Psi(u):\nabla u
      +\mathcal{PCL}\sigma_a(I_0-\Theta^4)-\mathcal{PL}(\sigma_a+\mathcal{L}_s\sigma_s)I_1\cdot u,  \\
  & \frac{\mathrm{St}}{\mathcal{C}} \partial_t I_0
    +\nabla \cdot  I_1=\mathcal{L}\sigma_a(\Theta^4-I_0),\\
  & \frac{\mathrm{St}}{\mathcal{C}} \partial_t I_1
     +\frac{1}{3}\nabla I_0=-\mathcal{L}(\sigma_a+\mathcal{L}_s\sigma_s)I_1,
\end{aligned}
\right.
\end{equation}
where we have dropped the subscript ``$*$" of the dimensionless quantities for the sake of simplicity.
And in the above dimensionless equations, we have used the following reduced dimensional parameters:
\begin{align*}
\mathrm{St} &= \frac{L_{\infty}}{t_{\infty}u_{\infty}},\!\!\!\!\!\!\!\!\!\!\!\!\!
  &&\mathrm{Strouhal}\;\mathrm{number},\;\;\;  \\
\mathrm{Ma}&= \frac{u_{\infty}}{\sqrt{\gamma R \Theta_{\infty}}},\!\!\!\!\!\!\!\!\!\!\!\!\!
 && \mathrm{Mach}\;\mathrm{number},\\
\mathrm{Re} &= \frac{L_{\infty}\rho_{\infty}u_{\infty}}{\mu_{\infty}},\!\!\!\!\!\!\!\!\!\!\!\!\!
 &&   \mathrm{Reynolds}\;\mathrm{number},\;\;\\
\mathrm{Pr}&= \frac{c_p\mu_{\infty}}{\kappa_{\infty}},\!\!\!\!\!\!\!\!\!\!\!\!\!
&&  \mathrm{Prandtl}\;\mathrm{number},\\
\mathcal{C}&= \frac{c}{u_{\infty}},\!\!\!\!\!\!\!\!\!\!\!\!\!
 && \mathrm{infrarelativisic}\;\mathrm{number},
\end{align*}
where
\begin{equation*}
   \gamma=\frac{c_p}{c_V}, \;\;\;c_p=R+c_V,\;\;
\end{equation*}
and
\begin{equation*}
  \mathcal{L}=L_{\infty}\sigma_{a,\infty},\quad \mathcal{L}_s=\frac{\sigma_{s,\infty}}{\sigma_{a,\infty}},\quad
  \mathcal{P}=\frac{a_r\Theta_{\infty}^4}{\rho_{\infty}e_{\infty}}.
\end{equation*}
From a physical point of view, $\mathcal{L}$, $\mathcal{L}\mathcal{L}_s$ and $\mathcal{P}$ measure  the strength of absorption, the strength of scattering, and the ratio of the radiative energy over the internal energy, respectively.

Our aim in this paper is to study the effects of the parameters $\mathrm{Ma}$, $\mathcal{P}$, $\mathcal{L}$, $\mathcal{L}_s$ and $\mathcal{C}$ on the radiation hydrodynamics model \eqref{A.1}.
Hence, we denote $\sqrt{\gamma}\mathrm{Ma}$ by $\mathrm{Ma}$ and take
$$\mathrm{St}=\mathrm{Re}=\sigma_a=\sigma_s=c_V=R=1\;\, \mathrm{and}\;\, \mathrm{Pr}=\gamma=2$$
in the above system to ignore the influence of these parameters and then obtain the dimensionless equations \eqref{NSFP1-Ma}.

\end{appendices}
\medskip
\medskip 
 \noindent
{\bf Acknowledgements:}
Li is supported by the National Natural Science Foundation of China (Grant No. 12331007), and the ``333 Project" of Jiangsu Province.
Zhang's work is supported by the National Natural Science Foundation of China  (Grant No.12501300), and the Basic Research Program of Jiangsu (Grant No. BK20250881).

\medskip
\noindent
{\bf  Data availability Statement:} Data sharing not applicable to this article as no datasets were generated or analysed during the current study.

\medskip
\noindent
{\bf  Conflict of interest statement:}
On behalf of all authors, the corresponding author states that there is no conflict of interest.

\bibliographystyle{plain}

\end{document}